\documentclass[10pt]{amsart}
\usepackage{latexsym,amsmath,amssymb,graphicx}
\usepackage[all,web]{xy}

\def\urlfont{\DeclareFontFamily{OT1}{cmtt}{\hyphenchar\font='057}
              \normalfont\ttfamily \hyphenpenalty=10000}

\DeclareFontFamily{OT1}{rsfs10}{}
\DeclareFontShape{OT1}{rsfs10}{m}{n}{ <-> rsfs10 }{}
\DeclareMathAlphabet{\mathscript}{OT1}{rsfs10}{m}{n}

\DeclareMathOperator{\im}{Im}       
\DeclareMathOperator{\id}{id}       
\DeclareMathOperator{\Pic}{Pic}     
\DeclareMathOperator{\Aut}{Aut}     
\DeclareMathOperator{\rk}{rk}       
\DeclareMathOperator{\Sing}{Sing}   
\DeclareMathOperator{\codim}{codim} 

\title[Geometric Transitions]{Geometric Transitions}

\author[Michele Rossi]{Michele Rossi}

\address{Dipartimento di Matematica, Universit\`a di Torino,
via Carlo Alberto 10, 10123 Torino} \email{michele.rossi@unito.it}

\thanks{Research partially supported by the
Italian PRIN project ``Geometria delle Variet\`{a} Algebriche"
(G.V.A.) and by the Dipartimento di Matematica dell'Universit\`{a}
di Torino (local research grant: ``Variet\`{a} Toriche,
transizioni geometriche e dualit\`{a} fisiche").}

\def\P{{\mathbb{P}}}

\def\p2{\mathbb{P}^2}
\def\p3{\mathbb{P}^3}
\def\p4{\mathbb{P}^4}

\def\su{\operatorname{SU}}

\def\rk{\operatorname{rk}}

\def\Z{\mathbb{Z}}

\def\C{\mathbb{C}}
\def\R{\mathbb{R}}

\def\Q{\mathbb{Q}}

\theoremstyle{plain}
\newtheorem{theorem}{Theorem}[section]
\newtheorem{proposition}[theorem]{Proposition}
\newtheorem{corollary}[theorem]{Corollary}
\newtheorem{conjecture}[theorem]{Conjecture}

\newtheorem*{conclusion}{Conclusion}

\theoremstyle{remark}
\newtheorem{remark}[theorem]{Remark}
\newtheorem{remarks}[theorem]{Remarks}
\newtheorem{example}[theorem]{Example}
\newtheorem{examples}[theorem]{Examples}
\newtheorem*{assumptions}{Assumptions}

\theoremstyle{definition}
\newtheorem{definition}[theorem]{Definition}

\newtheorem*{claim}{Claim}
\newtheorem*{step I}{Step I}
\newtheorem*{step II}{Step II}
\newtheorem*{step III}{Step III}
\newtheorem*{step IV}{Step IV}

\newcommand{\cy}{Ca\-la\-bi--Yau }
\newcommand{\ka}{K\"{a}hler }

\begin{document}

\begin{abstract}
The purpose of this paper is to give, on one hand, a mathematical
exposition of the main topological and geometrical properties of
geometric transitions, on the other hand, a quick outline of their
principal applications, both in mathematics and in physics.
\end{abstract}

\maketitle

\tableofcontents

\newpage
A geometric transition is a birational contraction followed by a
complex smoothing. This process connect two smooth, topologically
distinct, \cy threefolds. For this reason geometric transitions
attracted the interest of both mathematician and physicists.

\noindent From the mathematical point of view, the property of
changing topology candidates geometric transitions as the
3--dimensional analogous of analytic deformations between K3
surfaces. More precisely, K3 projective surfaces having different
sectional genus are linked by analytic deformations, showing that
their moduli space is actually connected. Analogously geometric
transitions may be the right way to give a notion of
``connectedness" to the ``moduli space" of \cy 3--folds. This is
essentially the famous Reid's fantasy \cite{Reid87} founded on
deep speculations due to H.~Clemens \cite{Clemens83}, R.~Friedman
\cite{Friedman86}, F.~Hirzebruch \cite{Hirzebruch87} and J.~Werner
\cite{Werner87}.

\noindent On the other hand, in physics, the same property
provides a mathematical tool to connect topologically distinct
compactifications to 4 dimensions of 10--dimensional type II
super--string theory vacua. This fact was firstly observed by
P.~Candelas, A.~M.~Dale, P.~S.~Green, T.~H\"{u}bsch,
C.~A.~L\"{u}tken and R.~ Schimmirgk in \cite{CDLS88},
\cite{Green-Hubsch88let}, \cite{Green-Hubsch88}, \cite{CGH89},
\cite{CGH90}. The physical interpretation of a geometric
transition connecting two topologically distinct string vacua was
given later, in 1995, by A.~Strominger \cite{Strominger95}, at
least in the case of a \emph{conifold} transition i.e. a geometric
transition whose associated birational contraction generates at
most ordinary double points. After this pivotal paper other
geometric transitions have been physically understood
\cite{BKK95}, \cite{KMP96}, \cite{BKKM97}.

\noindent For many geometric transitions, the induced change in
topology can be summarized by saying that \emph{a transition
increases complex moduli and decreases \ka moduli}. Since mirror
symmetry exchange complex and \ka moduli, it seemed natural to
conjecture the existence of a \emph{reverse transition} connecting
mirror partners of a couple of \cy 3-folds linked by a given
transition \cite{Morrison99}. Reverse transitions have been then
revealed useful tools for producing, at least conjecturally,
mirror constructions extending, via \emph{toric degenerations},
the Batyrev mirror symmetry between \cy 3--folds embedded in toric
varieties \cite{BC-FKvS98}, \cite{BC-FKvS00}, \cite{Batyrev04}.

\noindent In physics, geometric transitions have newly been in the
spot light as the geometric set up of recently conjectured
open/closed string dualities \cite{Gopakumar-Vafa99},
\cite{Ooguri-Vafa00}.

\vskip .3cm The present work is meant to give on one hand a
mathematical exposition of the main topological and geometrical
properties of a transition. This is the program of sections from 1
to 4: except for the latter, where some notion of deformation
theory in geometry is needed, these sections are devoted to
present a, as much as possible, self--contained treatment, for
graduate students and beginners. For this reason many well known
results or properties are developed in details like some example
(see \ref{l'esempio} and Example \ref{l'esempio bis}) or theorem
\ref{cambio omologico}. In particular the latter is intended to
give a complete account of the change in topology induced by a
conifold transition. Its content was already known twenty years
ago to H.~Clemens, and then to many other mathematicians and
physicists, but I was not able to find, in the literature, a
complete statement and a clear proof of all of the results
mentioned there. For this reason I preferred to rewrite here an
elementary proof requiring no more than basic facts in algebraic
topology and geometry.

\noindent On the other hand sections 5,6,7 give a quick outline of
some applications of geometric transitions both in mathematics and
in physics. Here the reader is clearly required to know basic
facts and definitions of these topics, although I tried to give
references of the original papers and, sometimes, of extended
surveys treating the mentioned subjects.

\vskip .3cm The paper is organized as follows.

\noindent Section 1 is devoted to give definition and examples of
geometric transitions. In particular the fundamental example of a
non--trivial conifold transition involving a quintic 3--fold in
$\P ^4$ is developed in detail.

\noindent Section 2 is a revised version of some of the
``topological considerations" given by H.~Clemens in
\cite{Clemens83}, which allow to locally think a conifold
transition as a surgery in topology (Proposition \ref{clemens
lemma}).

\noindent In section 3 the global change in topology induced by a
conifold transition is carefully studied, relying each other
homological invariants of all of the three poles of a conifold
transition (theorem \ref{cambio omologico}). This section ends up
with some similar considerations for more general geometric
transitions, essentially due to Y.~Namihawa and J.~Steenbrimk
\cite{Namikawa-Steenbrink95}.

\noindent Section 4 gives an outline of results and technics
needed to perform a (actually incomplete) classification of
geometric transition. Main results are here due to R.~Friedman,
M.~Gross and Y.~Namikawa.

\noindent The remaining sections are dedicated to describe some
fundamental applications of geometric transitions. Section 5
describes how geometric transitions are conjecturally employed, in
mathematics, to think of the \cy 3--folds moduli space as
``irreducible" and, in physics, to ``unify" type II super--string
compactified vacua. In section 6 a quick account of the role
played by geometric transitions in mirror symmetry is given,
starting from the key concept of reverse transition. In section 7
some further more recent applications are finally mentioned.

\vskip .4cm \noindent \textbf{Aknowledgements.} I would like to
especially thank A.~Collino and A.~Grassi for many suggestions and
stimulating discussions. In particular A.~Grassi introduced me to
this kind of problems. Thanks are also due to M.~Bill\'o and
I.~Pesando for many interesting discussions about string theory.

\section{Geometric Transitions: definition and the basic example}

\subsection{Calabi--Yau varieties}

\begin{definition}\label{cy-def}
Let $Y$ be a smooth, complex, projective variety with $\dim Y \geq
3$. $Y$ will be called a \emph{\cy variety} if
\begin{enumerate}
    \item $\bigwedge^{n} \Omega_{Y} =: \mathcal{K}_Y \cong \mathcal{O}_Y$
    \item $h^{p,0}(Y) = 0 \quad \forall 0<p<\dim Y$
\end{enumerate}
A 3--dimensional \cy variety will be also called a \emph{\cy
3--fold}.
\end{definition}

\begin{remarks}
\begin{enumerate}
    \item There are a lot of more or less equiv\-a\-lent def\-i\-ni\-tions
    of \cy varieties coming from:
    \begin{description}
        \item[differential geometry] the differential geometric concept  of a compact,
        \ka manifold admitting a \emph{Ricci flat metric} (Calabi conjecture and Yau
        theorem),
        \item[theoretical physics] the physical concept of a \ka , 3--dimensional
        complex, compact manifold admitting a flat,
        non--degenerate, holomorphic 3--form.
    \end{description}
    (see \cite{Joice2000} for a complete description of equivalences and
    implications).
    \item In the algebraic context, the given definition of \cy
    variety is the generalization of the following geometric
    objects
    \begin{description}
        \item[1--dimensional] smooth elliptic curves,
        \item[2--dimensional] smooth $K3$ surfaces.
    \end{description}

    \item With the dimensional bound $\dim Y \geq 3$, the given definition of \cy variety is equivalent
    to require that $Y$ is a \emph{\ka, compact, manifold whose holonomy group is a subgroup of $\su (\dim
    Y)$} (cfr. \cite{Joice2000}).

\end{enumerate}

\end{remarks}

\begin{examples}
\begin{enumerate}
    \item Smooth hypersurfaces of degree $n+1$ in $\P ^n$ (use Adjunction Formula and the Lefschetz Hyperplane Theorem).
    \item Smooth hypersurfaces (if exist!) of a weighted
    projective space $\P (q_0,\ldots,q_n)$ of degree
    $d=\sum_{i=0}^{n}q_i$.
    \item The general element of the anti--canonical system of a
    \emph{sufficiently good} 4--dimensional toric Fano variety (see \cite{Batyrev94}).
    \item Suitable complete intersections.... (iterate the previous
    examples).
    \item The double covering of $\P ^3$ ramified along a smooth
    surface of degree 8 in $\P ^3$ (octic double solid).
\end{enumerate}
\end{examples}

\subsection{Geometric transitions}

\begin{definition}(cfr. \cite{Morrison99}, \cite{Cox-Katz99}, \cite{GR02})
Let $Y$ be a \cy 3--fold and $\phi : Y\rightarrow \overline{Y}$ be
a \emph{birational contraction} onto a \emph{normal} variety. If
there exists a complex deformation (\emph{smoothing}) of
$\overline{Y}$ to a \cy 3--fold $\widetilde{Y}$, then the process
of going from $Y$ to $\widetilde{Y}$ is called a \emph{geometric
transition} (for short \emph{transition}) and denoted by
$T(Y,\overline{Y},\widetilde{Y})$ or by the diagram
\begin{equation*}
    \xymatrix@1{Y\ar@/_1pc/ @{.>}[rr]_T\ar[r]^{\phi}&
                \overline{Y}\ar@{<~>}[r]&\widetilde{Y}}\ .
\end{equation*}
A transition $T(Y,\overline{Y},\widetilde{Y})$ is called
\emph{trivial} if $\widetilde{Y}$ is a deformation of $Y$.
\end{definition}

\begin{remarks}
\begin{enumerate}
    \item \emph{Trivial transitions may occur}: e.g. consider Example
4.6 in \cite{Wilson92} where $\phi$ admits an elliptic scroll as
exceptional divisor and contracts it down to an elliptic curve
$C$.
    \item It is clearly possible to extend the transition process
    to any dimension $\geq 3$. Note that it is not possible to
    realize non--trivial transitions in dimension 1 (i.e. between
    elliptic curves).
    \item The transition process was firstly (locally) observed by
    H.~Clemens in the study of double solids $V$ admitting at worst
    nodal singularities \cite{Clemens83}: in his Lemma 1.11 he
    pointed out ``the relation of the resolution of the
    singularities of $V$ to the standard $S^3 \times D_3$ to $S^2 \times
    D_4$ surgery".
\end{enumerate}
\end{remarks}

\begin{definition}
A transition $T(Y,\overline{Y},\widetilde{Y})$ is called
\emph{conifold} if $\overline{Y}$ admits only \emph{ordinary
double points} (nodes) as singularities, i.e. singular points
whose tangent cones are singular hyperquadrics of rank $\dim X +1$
(precisely \emph{non--degenerate cones}).
\end{definition}

\subsection{The basic example: the conifold in $\P ^4$}\label{l'esempio}

The following example, given in \cite{GMS95}, shows that
\emph{non--trivial (conifold) transitions occur when $\dim X \geq
3$}.

Let $\overline{Y}\subset\P ^4$ be the singular hypersurface given
by the following equation
\begin{equation}\label{equazione}
    x_3 g(x_0,\ldots ,x_4) + x_4 h(x_0,\ldots ,x_4) = 0
\end{equation}
where $g$ and $h$ are generic homogeneous polynomials of degree 4.
$\overline{Y}$ is then the \emph{generic quintic 3--fold
containing the plane $\pi : x_3 = x_4 = 0$}. Then the singular
locus of $\overline{Y}$ is given by
\begin{equation}\label{singolarita}
    \Sing (\overline{Y}) = \{ [x]\in \P ^4 | x_3=x_4=g(x)=h(x)=0\}
\end{equation}

\begin{proposition}
$\Sing (\overline{Y})$ is composed by 16 nodes.
\end{proposition}

\begin{proof}
Let $p\in \Sing (\overline{Y})$. We have to write down the local
equation of $p$.

\noindent Assume $p=[1,0,0,0,0]$ and intersect $\overline{Y}$ with
the affine open subset of $\P ^4$
\begin{equation*}
    U_0 := \{ [x]\in\P ^4 | x_0\neq 0\}
\end{equation*}
Set $z_i:= x_i/x_0\ ,\ i=1,\ldots , 4$. Then $\overline{Y}\cap
U_0$ is described by the following affine equation
\begin{equation}\label{eq.affine}
    z_3 \widetilde{g}(z) + z_4 \widetilde{h}(z) = 0
\end{equation}
where $x_0^4\widetilde{g}=g$ and $x_0^4\widetilde{h}=h$. Besides
$p$ is the origin of $U_0$.

\noindent Since $g,h$ are generic we can assume that the
polynomial (holomorphic) maps $\widetilde{g},\widetilde{h}:\C
^4\rightarrow \C$ are submersive at the origin and we can find a
holomorphic chart $(U,z)$ centered in $p=0\in \C ^4$ and such that
\begin{equation}\label{eq.locale}
    \overline{U}:= \overline{Y}\cap U : z_3 z_1 + z_4 z_2 = 0
\end{equation}
Then $p$ is a node.
\end{proof}

\begin{proposition} [The resolution]\label{risoluzione}
$\Sing (\overline{Y})$ can be simultaneously resolved and the
resolution $\phi : Y\rightarrow \overline{Y}$ is a \emph{small
blow up} such that $Y$ is a smooth \cy 3--fold.
\end{proposition}

\begin{proof}
Blow up $\P ^4$ along the plane $\pi : x_3 = x_4 = 0$. We get a
birational morphism
\begin{equation*}
    \widehat{\phi}:\widehat{\P} ^4\longrightarrow \P ^4
\end{equation*}
whose exceptional divisor is a $\P ^1$--bundle over $\P ^2$. Let
$Y$ be the \emph{proper transform of} $\overline{Y}$ (i.e. the
closure in $\widehat{\P}^4$ of
$\widehat{\phi}^{-1}(\overline{Y}\setminus \pi)$). Since
$\widehat{\P} ^4$ is the hypersurface of bi--homogeneous equation
$y_0x_4 - y_1x_3 = 0$ in $\P ^4 (x) \times \P ^1 (y)$, then $Y$ is
the following complete intersection
\begin{eqnarray}\label{risol-equazioni}
  y_0x_4 - y_1x_3 &=& 0 \\
  \nonumber
  y_0 g(x) + y_1 h(x) &=& 0
\end{eqnarray}
and we get that
\begin{itemize}
    \item $Y$ is smooth,
    \item $\phi:= \widehat{\phi}_{|Y} : Y \longrightarrow
    \overline{Y}$ is an isomorphism outside of $\Sing
    (\overline{Y})$,
    \item $\forall p\in\Sing (\overline{Y})\quad \phi^{-1}(p) \cong \P
    ^1$.
\end{itemize}
Hence $\phi : Y \rightarrow\overline{Y}$ is a birational
resolution called \emph{small blow up} due to the dimension of its
exceptional locus ($1 < \dim Y -1 =2$).

\noindent To prove that $Y$ is \cy recall that $\widehat{\phi}$ is
a blow up, hence
\begin{equation*}
    K_{\widehat{\P} ^4} \equiv \widehat{\phi}^*(K_{\P ^4}) + (4-2-1)E \equiv
    - 5\widehat{\phi}^* (H)+ E
\end{equation*}
where $E$ is the exceptional divisor of $\widehat{\phi}$ and $H$
is the hyperplane of $\P ^4$. Then the Adjunction Formula gives
\begin{equation*}
    \mathcal{K}_{Y} \cong \mathcal{K}_{\widehat{\P} ^4}\otimes
    \mathcal{O}_{\widehat{\P}^4}(Y)\otimes \mathcal{O}_Y \cong \mathcal{O}_Y (E_{|Y}) \cong
    \mathcal{O}_Y
\end{equation*}
Moreover the Lefschetz Hyperplane Theorem and the K\"{u}nneth
Formula give
\begin{equation*}
  H^1 (Y,\C) \cong H^1 (\widehat{\P}^4,\C) \cong H^1(\P ^4\times\P ^1,\C)=0
\end{equation*}
hence $h^{1,0}(Y)=0$. On the other hand the Serre Duality theorem
allows to conclude that
\begin{equation*}
    H^2 (Y,\mathcal{O}_Y) \cong H^1 (Y,\mathcal{K}_Y) \cong H^1 (Y,\mathcal{O}_Y)
\end{equation*}
hence $h^{2,0}(Y) = h^{0,2} (Y) = h^{0,1}(Y) = h^{1,0}(Y)=0$.
\end{proof}

\begin{proposition} [The smoothing]
$\overline{Y}$ admits the obvious smoothing given by the generic
quintic 3--fold $\widetilde{Y}\subset \P ^4$. In particular
$\widetilde{Y}$ cannot be a deformation of $Y$ i.e. the conifold
transition $T(Y,\overline{Y},\widetilde{Y})$ is not trivial.
\end{proposition}

\begin{proof}
Apply again the Lefschetz Hyperplane Theorem and the K\"{u}nneth
Formula to get the following relations on the Betti numbers of
$\widetilde{Y}$ and $Y$
\begin{eqnarray}\label{betti_nmb}
\nonumber
  b_2 (\widetilde{Y})&=& b_2 (\P ^4) = 1 \\
  b_2 (Y) &=& b_2 (\P ^4 \times \P ^1) = 2
\end{eqnarray}
Therefore $\widetilde{Y}$ and $Y$ cannot be smooth fibers of the
same analytic family.
\end{proof}

\section{Local geometry and topology of a conifold
transition}\label{analisi locale}

The present section will be essentially devoted to explain the
basic argument given by  H.~Clemens in \cite{Clemens83}. As a consequence
we get that \emph{locally a conifold transition is described by a suitable surgery}.

\noindent In this section we will always
assume that $T(Y,\overline{Y},\widetilde{Y})$ is a \emph{conifold
transition} and $p$ is a point in $\Sing (\overline{Y})$, which means
that it is a \emph{node}.

\subsection{The local topology of a node}

\noindent Just like in the basic Example \ref{l'esempio}, we may
assume that there exists a local chart $(U,z)$ such that $p=0\in
U$. Denote $\overline{U}:=\overline{Y}\cap U$, which has local
equation in $U$ given by
\begin{equation}\label{equaz.locale}
    z_1z_3 + z_2z_4 = 0\ .
\end{equation}

\begin{proposition}\label{node-top}
Topologically $\overline{U}$ is a cone over $S^3 \times S^2$.
\end{proposition}

\begin{proof}
Change coordinates as follows
\begin{eqnarray}\label{cambio coordinate}
  w_1 &=& \frac{1}{2}(z_1 + z_3) \\
  \nonumber w_2 &=& \frac{i}{2}(-z_1 + z_3) \\
  \nonumber w_3 &=& \frac{1}{2}(z_2 + z_4) \\
  \nonumber w_4 &=& \frac{i}{2}(-z_2 + z_4)
\end{eqnarray}
to rewrite the local equation (\ref{equaz.locale}) as
\begin{equation*}
    \sum_{j=1}^{4} w_j^2 = 0 \ .
\end{equation*}
Decompose the latter in real and imaginary parts by setting $w_j =
u_j + iv_j$. Then $\overline{U}$ is described in $\R ^8 (u,v)$ by
the following two equations
\begin{eqnarray}\label{eq.i locali reali}
  \sum_{j=1}^{4} u_j^2 - \sum_{j=1}^{4} v_j^2  &=& 0 \\
  \sum_{j=1}^{4} u_j v_j  &=& 0 \ .
\end{eqnarray}
Fix now a real positive radius $\rho$ and consider the 7--sphere
\begin{equation*}
    S^{7}_{\rho}:=\{ (u,v)\in \R ^8 | \sum_{j=1}^{4} u_j^2 + \sum_{j=1}^{4} v_j^2 = \rho ^2 \}
\end{equation*}
Cut then $\overline{U}$ to get $\overline{U}_{\rho}:= \overline{U}
\cap S^7_{\rho}$. Topologically $\overline{U} =
\bigsqcup_{\rho\geq 0} \overline{U}_{\rho}$ and we get the claim
by proving that $\overline{U}_{\rho} \cong S^3 \times S^2$.

\noindent At this purpose, note that $\overline{U}_{\rho}$ is
described in $\R ^8$ by the following equations
\begin{eqnarray}\label{taglio sferico}
  \sum_{j=1}^{4} u_j^2 &=& \rho^2 -  \sum_{j=1}^{4} v_j^2\\
  \nonumber\sum_{j=1}^{4} v_j^2 &=& \frac{\rho^2}{2} \\
  \nonumber\sum_{j=1}^{4} u_j v_j   &=& 0
\end{eqnarray}
Then $\overline{U}_{\rho}$ can be fibred over the 3--sphere
$S^3_{\rho / \sqrt{2} }:=\{ v\in \R ^4 | \sum_{j=1}^{4} v_j^2 =
\rho ^2/2 \}$. Precisely the fiber over a point $v^o \in
S^3_{\rho/ \sqrt{2}}$ is given by
\begin{eqnarray*}
  \sum_{j=1}^{4} u_j^2 &=& \frac{\rho^2}{2} \\
  \sum_{j=1}^{4} v_j^o u_j  &=& 0
\end{eqnarray*}
which is a 2--sphere of radius $\rho /\sqrt{2}$.

\noindent The proof ends up by showing that the bundle
$\overline{U}_{\rho}$ is actually a product. This fact follows by
observing that  $\overline{U}_{\rho}$ is embedded in the tangent
bundle to the 3--sphere $S^3_{\rho / \sqrt{2} }\subset \R^4 (v)$.
In fact the latter is embedded in $\R^8(u,v)$ by the second and
third equations in (\ref{taglio sferico}). To conclude restrict to
$\overline{U}_{\rho}$ the well known trivialization
$T S^3_{\rho / \sqrt{2}}\cong S^3 \times \R ^3$.
\end{proof}

\subsection{Local geometry of the resolution}

To resolve the node recall Proposition \ref{risoluzione} of the
basic example. Precisely look at the \emph{proper transform}
$\widehat{U}$ of $\overline{U}$ in the blow up of the local chart
$(U,z)\cong \C ^4 (z)$ along the plane $z_3 = z_4 = 0$.

\noindent $\widehat{U}$ is then described in $\C ^4 \times \P ^1$ by the
following equations
\begin{eqnarray}\label{equazioni}
  y_0 z_4 - y_1 z_3 &=& 0 \\
  \nonumber y_0 z_1 + y_1 z_2 &=& 0
\end{eqnarray}

\begin{proposition}\label{risol-geom}
There is a diffeomorphism $\widehat{U}\cong \R ^4\times S^2$
\end{proposition}

\begin{proof}
Topologically it is not difficult to observe that $\widehat{U}$ is
an $\R ^4$--bundle over $\P ^1_{\C}$. In fact by splitting $z_j$
in real and imaginary parts, equations (\ref{equazioni}) give rise
to 4 linear equations in $\R ^8$ parameterized by $[y_0,y_1]\in \P
^1_{\C}$.

\noindent To construct the diffeomorphism introduce the
coordinates change given by (\ref{cambio coordinate}) and split
the new coordinates in real and imaginary parts: $w_j = u_j + i
v_j$. Equations (\ref{equazioni}) of $\widehat{U}$ can then be
rewritten in $\R ^8 (u,v)\times \P^1_{\C}$ in the following
matricial form
\begin{equation}
    \textbf{u}=A\left([y_0,y_1]\right)\textbf{v}
\end{equation}
where
\begin{equation}\label{matrice di Clemens}
    A\left([y_0,y_1]\right):= \left(%
\begin{array}{cccc}
  0 & |y_0|^2 - |y_1|^2 & 2 Im (\overline{y}_0 y_1) & 2 Re (\overline{y}_0 y_1) \\
  -|y_0|^2 + |y_1|^2 & 0 & -2 Re (\overline{y}_0 y_1) & -2 Im (\overline{y}_0 y_1) \\
  -2 Im (\overline{y}_0 y_1) & 2 Re (\overline{y}_0 y_1) & 0 & -|y_0|^2 + |y_1|^2 \\
  -2 Re (\overline{y}_0 y_1) & 2 Im (\overline{y}_0 y_1) & |y_0|^2 - |y_1|^2 & 0 \\
\end{array}%
\right)
\end{equation}
We will refer to the matrix $A$ as the \emph{Clemens' matrix}: in
fact it is the same matrix appearing in formula (1.18) of
\cite{Clemens83}. For any $[y]\in \P ^1_{\C}$, one can easily
check that $A[y]\in SO (4)$ and moreover it is antisymmetric i.e.
$^{t}A[y] + A[y] = 0$.

\noindent A diffeomorphism $\Phi:\widehat{U}\cong \R ^4 \times
\P^1_{\C}$ is then given by
\begin{equation}\label{diffeomorfismo}
\begin{array}{cccc}
  \Phi^{-1}: & \R ^4 \times \P^1_{\C} & \longrightarrow & \widehat{U}\subset \R^4 \times \R^4 \times \P^1_{\C} \\
   & (v,[y]) & \mapsto & (A[y]v,v,[y]) \\
\end{array}
\end{equation}
The proof ends up by the usual identification $\P^1_{\C}\cong
S^2$.
\end{proof}

\begin{remark}Just like in the basic Example \ref{l'esempio}, the
restriction of the blow up of $U=\C ^4$ along the plane
$z_3=z_4=0$ gives rise to a birational map
\begin{equation*}
\xymatrix@1{\varphi: \widehat{U}\ar[r]&\overline{U}}
\end{equation*}
which is a \emph{small blow up}. Precisely $\varphi$ is biregular
over the complement of the origin in $\overline{U}$ and
$\varphi^{-1}(0) = \P^1_{\C}$. Then it induces a diffeomorphism
\begin{equation*}
\xymatrix@1{ \widehat{U}\setminus\varphi^{-1}(0)
\ar[r]^{\varphi}_{\cong} & \overline{U}\setminus\{0\} }
\end{equation*}
Recalling Proposition \ref{node-top},
$\overline{U}\setminus\{0\}\cong (\R ^4)\setminus\{0\} \times
\P^1_{\C}$ and it is natural to ask what is the relation between
$\varphi$ and $\Phi$. Thanks to the Clemens matrix's properties we
get that
\begin{equation}\label{estensione}
    \Phi|_{\widehat{U}\setminus\varphi^{-1}(0)} = \varphi |_{\widehat{U}\setminus\varphi^{-1}(0)}
\end{equation}
and $\Phi$ \emph{is an extension of $\varphi$ over the exceptional
fibre} i.e. the following commutative diagram holds
\begin{equation*}
\xymatrix{\widehat{U} \ar[r]^\Phi & \R^4 \times S^2 \\
\widehat{U}\setminus\varphi^{-1}(0) \ar@{^{(}->}[u]
\ar[r]^{\varphi} & \overline{U}\setminus\{0\} \ar@{^{(}->}[u]}
\end{equation*}
To prove this fact it suffices to check that $(u,v)=(Av,v)$
satisfies the real equations (\ref{eq.i locali reali}) of
$\overline{U}$, for any $v\neq 0$. In fact
\begin{equation*}
|u|^2 - |v|^2 = ^{t}v\ ^{t}A\ A\  v -\  ^{t}v\ v = 0
\end{equation*}
since $A$ is orthogonal. On the other hand
\begin{equation*}
    \sum_{j=1}^{4} u_j v_j = ^{t}v\ ^{t}A\ v = -\ ^{t}v\ A\ v = 0
\end{equation*}
since $A$ is antisymmetric and it induces an alternating bilinear
form.
\end{remark}

\begin{proposition}
$\widehat{U}$ can be identified with the total space of the rank 2
holomorphic vector bundle
$\mathcal{O}_{\P^1}(-1)\oplus\mathcal{O}_{\P^1}(-1)$ over the
exceptional fibre $\P^1_{\C}=\varphi^{-1}(0)$. In particular
$\widehat{U}$ admits a natural complex structure.
\end{proposition}

\begin{proof}
Since $\widehat{U}$ is the proper transform of $\overline{U}$ in
the small blow up of $U$, it can be identified with the total
space of the normal bundle of the exceptional fibre $\mathcal{N}
_{\widehat{U}|\mathbb{P}^{1}}$. The latter is a holomorphic vector
bundle of rank 2 over the exceptional fibre $\mathbb{P}^{1}_{\C}$.
By the Grothendieck theorem it splits as follows
\begin{equation*}
\mathcal{N}_{\widehat{U}|\mathbb{P}^{1}}\cong
\mathcal{O}_{\mathbb{P}^{1}}(d_{1})\oplus
\mathcal{O}_{\mathbb{P}^{1}}(d_{2})\ .
\end{equation*}
Choose two local charts on $S^{2}\cong \mathbb{P}^{1}_{\C}$ around
the north and the south poles respectively. Let $\tau
:=y_{0}/y_{1}$ and $\sigma :=y_{1}/y_{0}$ be the associated local
coordinates. Lifting these charts to
$\mathcal{O}_{\mathbb{P}^{1}}(d_{1})\oplus
\mathcal{O}_{\mathbb{P}^{1}}(d_{2})$ means that we can choose two
local parameterizations
\begin{equation*}
(\tau ;t_{1},t_{2})\quad ,\quad (\sigma ;s_{1},s_{2})
\end{equation*}
patching along the fibre over the fixed point $(y_{0}:y_1)=(\tau
:1)=(1:\sigma )$ as follows
\begin{equation*}
s_{i}=\tau ^{-d_{i}}t_{i}
\end{equation*}
where $\tau ^{-d_{i}}$ represents the transition function in $GL(1,\mathbb{C})=%
\mathbb{C}^{*}$.

\noindent Equations (\ref{equazioni}) of $\widehat{U}$ allow us to
set
\begin{equation*}
t_{1}=z_{1}\ ,\ t_{2}=z_{4}\ ;\ s_{1}=-z_{2}\ ,\ s_{2}=z_{3}
\end{equation*}
Then
\begin{eqnarray*}
s_{1} &=&-z_{2}=\frac{y_{0}}{y_{1}}z_{1}=\tau t_{1} \\
s_{2} &=&z_{3}=\frac{y_{0}}{y_{1}}z_{4}=\tau t_{2}
\end{eqnarray*}
and we get that $d_{1}=d_{2}=-1$.
\end{proof}

\subsection{Local geometry of the smoothing}

Recalling the real equations (\ref{eq.i locali reali}) of
$\overline{U}$, a local smoothing of the node is given by the
1--parameter family $f:\mathcal{U}\rightarrow\R$ where
\begin{equation}\label{lisciamento}
    U_t:=f^{-1}(t):\left\{\begin{array}{c}
      \sum_{j=1}^{4} u_j^2 - \sum_{j=1}^{4} v_j^2  = t \\
      \sum_{j=1}^{4} u_j v_j  = 0 \\
    \end{array}\right.
\end{equation}
Let $\widetilde{U}:=U_{t_0}$ for some $t_0 \in \R, t_0>0$.

\begin{proposition}\label{smooth-geom}
$\widetilde{U}$ is diffeomorphic to the cotangent bundle $T^* S^3$
of the 3--sphere. In particular $\widetilde{U}\cong
S^3\times\R^3$.
\end{proposition}

\begin{proof}
$T^* S^3$ can be embedded in $\R^8(q,p)$ by the standard equations
\begin{eqnarray*}
  \sum_{j=1}^{4} q_j^2 &=& 1 \\
  \sum _{j=1}^{4} q_j p_j &=& 0 \ .
\end{eqnarray*}
The diffeomorphism $\Psi : \widetilde{U}\cong T^* S^3$ is then
defined by setting
\begin{eqnarray*}
  q_j &=& \frac{u_j}{\sqrt{t_0 + \sum_{j} v_j^2}} \\
  p_j &=& v_j
\end{eqnarray*}
The proof concludes by applying the standard trivialization $T^*
S^3\cong S^3\times\R^3$.
\end{proof}

\begin{remark}\label{tildeS}
The vanishing cycle of the smoothing $f:\mathcal{U}\rightarrow\R$
is given by the family of embedded 3--spheres
$\mathcal{S}\rightarrow\R$ defined by
\begin{equation}\label{ciclo evenescente}
    S_t:=\left\{
\begin{array}{cc}
  |u|^2 - t = v_1 = \ldots = v_4 = 0 & \text{if}\ t\geq 0 \\
  |v|^2 - t = u_1 = \ldots = u_4 = 0 & \text{if}\ t\leq 0 \\
\end{array}%
\right.
\end{equation}
Clearly $S_0 = \{0\}\subset\overline{U}$. Define $\widetilde{S}:=
S_{t_0}$. Recalling the diffeomorphism $\Psi$ of the previous
proposition we get that $\Psi (\widetilde{S})$ \emph{is the
0--section of the cotangent bundle} $T^* S^3$.
\end{remark}

\begin{definition}
Let $L$ be a submanifold of a given symplectic manifold
$(M,\omega)$. $L$ is called \emph{lagrangian} if
\begin{enumerate}
    \item $2 \dim_{\R}L = \dim_{\R}M$
    \item $\forall p \in L , \forall X,Y \in T_p M ,\quad  \omega_p
    (X,Y)=0$ \ .
\end{enumerate}
\end{definition}

\begin{example}\label{0-sezione}
The cotangent bundle $T^* M$ of a given manifold $M$ admits the
\emph{ canonical}  symplectic structure given by
$\omega:=d\vartheta$, where $\vartheta$ is the Liouville 1--form.
The 0--section of $T^* M$ is a lagrangian submanifold with respect
to the canonical symplectic structure.
\end{example}

\begin{proposition}
$\widetilde{U}$ admits a natural symplectic structure and the
vanishing cycle $\widetilde{S}$ is a lagrangian submanifold.
\end{proposition}

\begin{proof}
Let $\omega$ be the canonical symplectic structure on $T^* S^3$.
Then $\Psi^*(\omega)$ gives the natural symplectic structure to
$\widetilde{U}$. By remark \ref{tildeS} and Example
\ref{0-sezione} we get that
\begin{equation*}
    \Psi^*(\omega)|_{\widetilde{S}} = \omega |_{S^3} = 0 \ .
\end{equation*}
\end{proof}

\subsection{Local topology of a conifold transition}

\begin{proposition}[\cite{Clemens83}, Lemma 1.11]\label{clemens lemma}
Let $D_n\subset \R^n$ be the closed unit ball and consider
\begin{itemize}
    \item $S^3\times D_3\subset S^3\times \R^3 \overset{\Psi^{-1}}{\cong} \widetilde{U}$
    \item $D_4\times S^2\subset \R^4\times S^2 \overset{\Phi^{-1}}{\cong} \widehat{U}$
\end{itemize}
Then $\widetilde{D}:= \Psi^{-1}(S^3\times D_3)$ and $\widehat{D}:=
\Phi^{-1}(D_4\times S^2)$ are compact tubular neighborhoods of the
vanishing cycle $\widetilde{S}\subset\widetilde{U}$ and of the
exceptional cycle $\P^1_{\C}\subset\widehat{U}$, respectively.

\noindent Consider the standard diffeomorphism
\begin{equation*}
    \begin{array}{cccc}
      \alpha' : & (\R^4 \setminus \{0\})\times S^2 & \overset{\cong}{\longrightarrow} & S^3 \times (\R^3 \setminus \{0\}) \\
       & (u,v) & \mapsto &  ( \frac{u}{|u|}, |u|v ) \\
    \end{array}
\end{equation*}
and restrict it to $D_4 \times S^2$. Since
\begin{equation*}
    \partial(D_4 \times S^2)= S^3 \times S^2 = \partial(S^3 \times D_3)
\end{equation*}
observe that $\alpha'|_{\partial(D_4 \times S^2)}= \id |_{S^3
\times S^2}$. Hence $\alpha'$ induces a standard surgery from
$\R^4\times S^2$ to $S^3\times \R^3$.

\noindent Then $\widetilde{U}$ can be obtained from $\widehat{U}$
by removing $\widehat{D}$ and pasting in $\widetilde{D}$, by means
of the diffeomorphism $\alpha:= \Psi^{-1}\circ\alpha'\circ\Phi$.
\end{proposition}

\begin{proof}
The situation is described by the following commutative diagram
\begin{equation*}
    \xymatrix{\widehat{U}\setminus \P^1_{\C}
    \ar[r]^{\alpha}_{\cong} \ar[d]^{\Phi=\varphi}_{\cong}&
    \widetilde{U}\setminus\widetilde{S}\ar[d]^{\Psi}_{\cong}\\
    (\R^4 \setminus \{0\})\times S^2 \ar[r]^{\alpha'}_{\cong} &
    S^3 \times (\R^3 \setminus \{0\})}
\end{equation*}
which implies that $\alpha$ induces a diffeomorphism from
$\partial(\widetilde{D})$ to $\partial(\widehat{D})$. The claim
follows immediately.
\end{proof}

\section{Global geometry and topology of a conifold transition}

Let $T(Y,\overline{Y},\widetilde{Y})$ be a \emph{conifold}
transition. Then, by definition and the local analysis of the
previous section we know that:
\begin{itemize}
    \item $\Sing (\overline{Y})= \{ p_1,\ldots , p_N\}$ where
    $p_i$ is a node;
    \item there exists a simultaneous resolution
    $\phi:Y\rightarrow\overline{Y}$ which is a birational morphism
    contracting $N$ rational curves $E_1, \ldots, E_N$;
    \item $\widetilde{Y}$ admits $N$ vanishing cycles $S_1, \ldots,
    S_N$ which are 3-spheres.
\end{itemize}
Two natural questions then arise:
\begin{enumerate}
    \item Are the homology classes $[E_1],\ldots ,[E_N]\in H_2(Y,
    \Z)$ linearly independent? Which is: are the exceptional
    curves of $\phi$ homologically independent?
    \item Same question about $[S_1],\ldots ,[S_N]\in H_3(\widetilde{Y},
    \Z)$, i.e. are the vanishing cycles homologically independent?
\end{enumerate}
The answer is \emph{no} to both questions!

\begin{example}\label{l'esempio bis}
Consider the example given in \ref{l'esempio} of the conifold in
$\P^4$. Then $\overline{Y}=\{ x_3g+x_4h=0\}\subset\P^4$, $N=16$,
the resolution $Y$ contains 16 exceptional rational curves and the
smoothing $\widetilde{Y}$ contains 16 vanishing spheres.

\noindent For question (1) notice that if $[E_1],\ldots ,[E_{16}]$
would be independent then we would have
\begin{equation*}
    b_2(Y)=b_2(\widetilde{Y})+16
\end{equation*}
which is clearly contradicting (\ref{betti_nmb}).

\noindent On the other hand, for question (2) let us compare
$b_3(Y)$ and $b_3(\widetilde{Y})$.

\begin{claim}
$b_3(Y)=174$, $b_3(\widetilde{Y})=204$; then
$b_3(\widetilde{Y})-b_3(Y)=30$.
\end{claim}

\begin{proof}
In physics literature, this proof is often realized by invoking
the local smoothness of the complex moduli space of a \cy 3--fold
$Y$, hence the Bogomolov--Tian--Todorov theorem (see
\cite{Bogomolov78}, \cite{Tian87}, \cite{Todorov89}, \cite{Ran92};
see also the following \ref{cy-moduli}). Then it is well defined a
tangent space, to such a moduli space, canonically identified with
$H^1(\mathcal{T}_Y)$, via the Kodaira--Spencer map. The \cy
condition gives then
\begin{equation*}
    b_3(Y) = 2 + 2 h^{2,1}(Y) = 2 + 2h^1(\mathcal{T}_Y)
    \ .
\end{equation*}
The statement follows, for both $Y$ and $\widetilde{Y}$, by
counting their moduli (see \cite{GMS95}).

\noindent Actually proving the claim do not need local smoothness
of moduli spaces, which is a very deeper concept. In the following
we present a more (standard) elementary proof. Although
computationally more intricate than the previous one, the
following method has the advantage to apply to more general
situations: in fact it is not easy to count moduli of a general
\cy 3--fold, even in the case of a complete intersection.

\noindent Let start to consider $\widetilde{Y}$ which is the
easiest case of a projective hypersurface. In this case there are
many methods to compute $h^{2,1}(\widetilde{Y})$: e.g. it is
possible to compute directly $h^1(\mathcal{T}_Y)$ by Poincar\'e
residues (see \cite{Griffiths69}) and to end up by using \cy
condition. Here is the most elementary procedure to compute
$h^1(\mathcal{T}_Y)$.

\noindent Since $\mathcal{N}_{\widetilde{Y}\mid \P^4}\cong
\mathcal{O}_{\P^4}(5)\otimes\mathcal{O}_{\widetilde{Y}}=:\mathcal{O}_{\widetilde{Y}}(5)$,
the tangent sheaf exact sequence gives
\begin{equation*}
    \xymatrix@1{0\ar[r]&\mathcal{T}_{\widetilde{Y}}\ar[r]&
                \mathcal{T}_{\P^4}\otimes\mathcal{O}_{\widetilde{Y}}\ar[r]&
                \mathcal{O}_{\widetilde{Y}}(5)\ar[r]&0}
\end{equation*}
and the associated cohomology long exact sequence starts as
follows
\begin{equation}\label{coom-succ-fs-tg}
    0\rightarrow H^0 \left (\mathcal{T}_{\widetilde{Y}}\right )\rightarrow
                H^0\left (\mathcal{T}_{\P^4}\otimes\mathcal{O}_{\widetilde{Y}}\right )\rightarrow
                H^0\left (\mathcal{O}_{\widetilde{Y}}(5)\right )\rightarrow
                H^1\left (\mathcal{T}_{\widetilde{Y}}\right )\rightarrow
                H^1\left (\mathcal{T}_{\P^4}\otimes\mathcal{O}_{\widetilde{Y}}\right
                )\rightarrow
\end{equation}
All needed information can then be deduced by the cohomology
associated with the Euler exact sequence
\begin{equation}\label{Euler succ}
    \xymatrix@1{0\ar[r]&\mathcal{O}_{\widetilde{Y}}\ar[r]&
                \mathcal{O}_{\P^4}(1)^{\oplus 5}\ar[r]&
                \mathcal{T}_{\P^4}\otimes\mathcal{O}_{\widetilde{Y}}\ar[r]&0}
\end{equation}
and with the following tensor product, by
$\mathcal{O}_{\widetilde{Y}}(5)$, of the structure sheaf exact
sequence of $\widetilde{Y}\subset\P^4$
\begin{equation}\label{succ-fs-str}
    \xymatrix@1{0\ar[r]&\mathcal{O}_{\widetilde{Y}}\ar[r]&
                \mathcal{O}_{\P^4}(5)\ar[r]&
                \mathcal{O}_{\widetilde{Y}}(5)\ar[r]&0} \ .
\end{equation}
In fact (\ref{Euler succ}) gives
\begin{eqnarray*}
    \xymatrix@1{0\ar[r]&\C\ar[r]&H^0\left (\mathcal{O}_{\P^4}(1)\right )^{\oplus 5}\ar[r]&
              H^0\left(\mathcal{T}_{\P^4}\otimes\mathcal{O}_{\widetilde{Y}}\right)\ar[r]&
              H^1\left(\mathcal{O}_{\widetilde{Y}}\right)\ar[r]&}&&\\
    \nonumber
    \xymatrix@1{\ar[r]&H^1\left (\mathcal{O}_{\P^4}(1)\right )^{\oplus 5} \ar[r]&
              H^1\left(\mathcal{T}_{\P^4}\otimes\mathcal{O}_{\widetilde{Y}}\right)\ar[r]&
              H^2\left(\mathcal{O}_{\widetilde{Y}}\right)\ar[r]&\cdots}&&
\end{eqnarray*}
Bott formulas
\begin{equation}\label{bott}
    h^q\left(\Omega^p_{\P^n}(a)\right)=\left\{\begin{array}{cc}
      \binom{a+n-p}{a}\binom{a-1}{p} & \text{for $q=0$, $0\leq p\leq n$ and $a>p$,} \\
      1 & \text{for $0\leq p = q\leq n$ and $a=0$,} \\
      \binom{-a+p}{-a}\binom{-a-1}{n-p} & \text{for $q=n$, $0\leq p\leq n$ and $a<p-n$,} \\
      0 & \text{otherwise} \\
    \end{array}\right.
\end{equation}
and the \cy condition
$h^1(\mathcal{O}_{\widetilde{Y}})=h^2(\mathcal{O}_{\widetilde{Y}})=0$,
allow then to conclude that
\begin{eqnarray*}
  h^0 \left(\mathcal{T}_{\P^4}\otimes\mathcal{O}_{\widetilde{Y}}\right) &=& 25 -1 = 24 \\
  \nonumber
  h^1 \left(\mathcal{T}_{\P^4}\otimes\mathcal{O}_{\widetilde{Y}}\right) &=& 0
\end{eqnarray*}
On the other hand the cohomology of (\ref{succ-fs-str}) gives
\begin{equation*}
    \xymatrix@1{0\ar[r]&\C\ar[r]&H^0 \left(\mathcal{O}_{\P^4}(5)\right)\ar[r]&
                H^0 \left(\mathcal{O}_{\widetilde{Y}}(5)\right)\ar[r]&
                H^1 \left(\mathcal{O}_{\widetilde{Y}}\right)\ar[r]& \cdots }
\end{equation*}
Again Bott formulas (\ref{bott}) and \cy condition imply that
\begin{equation*}
    h^0\left(\mathcal{O}_{\widetilde{Y}}(5)\right) = 126 -1 =125 \ .
\end{equation*}
Since $h^0 \left(\mathcal{T}_{\widetilde{Y}}\right)=
h^0\left(\Omega^2_{\widetilde{Y}}\right)=0$, the sequence
(\ref{coom-succ-fs-tg}) gives
\begin{equation*}
    h^1 \left(\mathcal{T}_{\widetilde{Y}}\right ) = 125 - 24 = 101\ .
\end{equation*}
The previous argument do not apply to the resolution $Y$, since it
is the complete intersection given by the bi--homogeneous
equations (\ref{risol-equazioni}) in $\P^1 \times \P^4 =: \P$. In
this case there is no more an Euler sequence like (\ref{Euler
succ}), then it is better to directly compute $h^1(\Omega^2_Y)$.
At this purpose dualize the tangent sheaf sequence to get
\begin{equation}\label{diff-succ}
    \xymatrix@1{0\ar[r]&\mathcal{N}^*_{Y\mid\P}\ar[r]&
                \Omega_{\P}\otimes\mathcal{O}_Y\ar[r]&
                \Omega_Y\ar[r]&0}
\end{equation}
where $\mathcal{N}^*_{Y\mid\P}:=\mathcal{H}om
\left(\mathcal{N}_{Y\mid\P},\mathcal{O}_Y\right)=\mathcal{I}_Y/\mathcal{I}_Y^2$,
being $\mathcal{I}_Y$ the ideal sheaf of $Y\subset\P$. Then
\begin{equation*}
    \mathcal{N}^*_{Y\mid\P}\cong
    \left[\mathcal{O}_{\P}(-1,-1)\oplus\mathcal{O}_{\P}(-1,-4)\right]\otimes\mathcal{O}_Y
    =:\mathcal{O}_Y(-1,-1)\oplus\mathcal{O}_Y(-1,-4) \ .
\end{equation*}
Since $Y$ is \cy, its canonical sheaf is trivial and the fourth
exterior power of (\ref{diff-succ}) gives the following exact
sequence
\begin{equation*}
    \xymatrix@1{0\ar[r]&\mathcal{O}_Y(-2,-5)\otimes \Omega^2_Y\ar[r]&
                \Omega_{\P}^4\otimes\mathcal{O}_Y\ar[r]&
                \mathcal{O}_Y(-1,-1)\oplus\mathcal{O}_Y(-1,-4)\ar[r]&0}
\end{equation*}
This sequence, tensored by $\mathcal{O}_{\P}(2,5)$, gives then
rise to the following one
\begin{equation}\label{successione}
    \xymatrix@1{0\ar[r]&\Omega^2_Y\ar[r]&
                \Omega_{\P}^4(2,5)\otimes\mathcal{O}_Y\ar[r]&
                \mathcal{O}_Y(1,4)\oplus\mathcal{O}_Y(1,1)\ar[r]&0}
\end{equation}
from which it is possible to compute $h^1(\Omega^2_Y)$ by passing
to the associated long exact sequence in cohomology. In fact,
recalling the \cy condition for $Y$, it follows that
\begin{eqnarray}\label{coom-successione}
    &&\xymatrix@1{0\ar[r]&H^0 \left(\Omega_{\P}^4(2,5)\otimes\mathcal{O}_Y\right)\ar[r]&
                H^0\left(\mathcal{O}_Y(1,4)\right)\oplus
                H^0\left(\mathcal{O}_Y(1,1)\right)\ar[r]&}\\
    \nonumber
    &&\xymatrix@1{\ar[r]& H^1 \left(\Omega^2_Y\right)\ar[r]&
                H^1\left(\Omega_{\P}^4(2,5)\otimes\mathcal{O}_Y\right)\ar[r]&\cdots}
\end{eqnarray}
All needed information can then be obtained by suitable twists of
the following structure sheaves exact sequences of
$Y\subset\widehat{\P}\subset\P$
\begin{equation}\label{succI-fs-str}
    \xymatrix@1{0\ar[r]&\mathcal{O}_{\P}(-1,-1)\ar[r]&\mathcal{O}_{\P}\ar[r]&
                \mathcal{O}_{\widehat{\P}}\ar[r]&0}
\end{equation}
\begin{equation}\label{succII-fs-str}
    \xymatrix@1{0\ar[r]&\mathcal{O}_{\widehat{\P}}(-1,-4)\ar[r]&
                \mathcal{O}_{\widehat{\P}}\ar[r]& \mathcal{O}_Y\ar[r]&0}
\end{equation}
where $\widehat{\P}$ is the blow up of $\P^4$ along the plane
$x_3=x_4=0$, whose equation in $\P$ is the former in
(\ref{risol-equazioni}), and $\mathcal{O}_{\widehat{\P}}(-1,-4):=
\mathcal{O}_{\P}(-1,-4)\otimes \mathcal{O}_{\widehat{\P}}$\ .

\noindent In fact the tensor product of (\ref{succI-fs-str}) and
(\ref{succII-fs-str}) by $\Omega_{\P}^4(2,5)$ gives
\begin{equation}\label{Omega-twistI}
    \xymatrix@1{0\ar[r]&\Omega_{\P}^4(1,4)\ar[r]&\Omega_{\P}^4(2,5)\ar[r]&
                  \Omega_{\P}^4(2,5)\otimes\mathcal{O}_{\widehat{\P}}\ar[r]&0}
\end{equation}
\begin{equation}\label{Omega-twistII}
    \xymatrix@1{0\ar[r]&\Omega_{\P}^4(1,1)\otimes\mathcal{O}_{\widehat{\P}}\ar[r]&
                  \Omega_{\P}^4(2,5)\otimes\mathcal{O}_{\widehat{\P}}\ar[r]&
                  \Omega_{\P}^4(2,5)\otimes\mathcal{O}_Y\ar[r]&0}
\end{equation}
The following K\"unneth formulas
\begin{equation}\label{kunneth}
    h^v\left(\Omega^u_{\P}(a,b)\right) = \bigoplus_{\text{
      $\begin{array}{c}
        p+r = u \\
        q+s = v \\
      \end{array}$}}
    \left[h^q\left(\Omega^p_{\P^1}(a)\right)\cdot
                       h^s\left(\Omega^r_{\P^4}(b)\right)\right]
\end{equation}
and (\ref{bott}) applied to the cohomology long exact sequence of
(\ref{Omega-twistI}) give
\begin{eqnarray*}
  h^0 \left(\Omega_{\P}^4(2,5)\otimes\mathcal{O}_{\widehat{\P}}\right) &=&
  h^0 \left(\Omega_{\P}^4(2,5)\right) - h^0\left(\Omega_{\P}^4(1,4)\right) = 27\\
  h^1 \left(\Omega_{\P}^4(2,5)\otimes\mathcal{O}_{\widehat{\P}}\right)&=& 0
\end{eqnarray*}
Moreover the tensor product of (\ref{succI-fs-str}) by
$\Omega_{\P}^4(1,1)$ gives
\begin{equation}\label{Omega-twistIII}
    \xymatrix@1{0\ar[r]&\Omega_{\P}^4\ar[r]&\Omega_{\P}^4(1,1)\ar[r]&
                  \Omega_{\P}^4(1,1)\otimes\mathcal{O}_{\widehat{\P}}\ar[r]&0}
\end{equation}
whose cohomology attains the following results
\begin{equation*}
    h^0\left(\Omega_{\P}^4(1,1)\otimes\mathcal{O}_{\widehat{\P}}\right)=
    h^1\left(\Omega_{\P}^4(1,1)\otimes\mathcal{O}_{\widehat{\P}}\right)=0
    \ .
\end{equation*}
Therefore the cohomology of (\ref{Omega-twistII}) allows to
conclude that
\begin{eqnarray}\label{risultatiI}
  h^0 \left(\Omega_{\P}^4(2,5)\otimes\mathcal{O}_Y\right) &=&
  h^0 \left(\Omega_{\P}^4(2,5)\otimes\mathcal{O}_{\widehat{\P}}\right)= 27\\
\nonumber
  h^1 \left(\Omega_{\P}^4(2,5)\otimes\mathcal{O}_Y\right)&=&
  h^1 \left(\Omega_{\P}^4(2,5)\otimes\mathcal{O}_{\widehat{\P}}\right)=0
\end{eqnarray}
To compute $h^0(\mathcal{O}_Y(1,4))$, consider the tensor product
of (\ref{succI-fs-str}) and (\ref{succII-fs-str}) by
$\mathcal{O}_Y(1,4)$:
\begin{equation*}
    \xymatrix@1{0\ar[r]&\mathcal{O}_{\P}(0,3)\ar[r]&\mathcal{O}_{\P}(1,4)\ar[r]&
                \mathcal{O}_{\widehat{\P}}(1,4)\ar[r]&0}
\end{equation*}
\begin{equation*}
    \xymatrix@1{0\ar[r]&\mathcal{O}_{\widehat{\P}}\ar[r]&
                \mathcal{O}_{\widehat{\P}}(1,4)\ar[r]& \mathcal{O}_Y(1,4)\ar[r]&0}
\end{equation*}
Again formulas (\ref{kunneth}) and (\ref{bott}) applied to the
cohomology of the first sequence give
\begin{equation*}
    h^0\left(\mathcal{O}_{\widehat{\P}}(1,4)\right)=
    h^0\left(\mathcal{O}_{\P}(1,4)\right)-
    h^0\left(\mathcal{O}_{\P}(0,3)\right)= 140 -35 = 105 \ .
\end{equation*}
The cohomology of the second sequence allows to conclude
\begin{equation}\label{risultatiII}
    h^0(\mathcal{O}_Y(1,4)) = 104 \ .
\end{equation}
Analogously for $h^0(\mathcal{O}_Y(1,1))$ one has
\begin{equation*}
    \xymatrix@1{0\ar[r]&\mathcal{O}_{\P}\ar[r]&\mathcal{O}_{\P}(1,1)\ar[r]&
                \mathcal{O}_{\widehat{\P}}(1,1)\ar[r]&0}
\end{equation*}
\begin{equation*}
    \xymatrix@1{0\ar[r]&\mathcal{O}_{\widehat{\P}}(0,-3)\ar[r]&
                \mathcal{O}_{\widehat{\P}}(1,1)\ar[r]&
                \mathcal{O}_Y(1,1)\ar[r]&0}\ .
\end{equation*}
Then
\begin{equation*}
    h^0\left(\mathcal{O}_{\widehat{\P}}(1,1)\right)=
    h^0\left(\mathcal{O}_{\P}(1,1)\right)-
    h^0\left(\mathcal{O}_{\P}\right)= 10 -1 = 9
\end{equation*}
and finally
\begin{equation}\label{risultatiIII}
    h^0(\mathcal{O}_Y(1,1)) = 9 \ .
\end{equation}
Therefore, recalling (\ref{coom-successione}), results
(\ref{risultatiI}), (\ref{risultatiII}) and (\ref{risultatiIII})
end up the proof giving
\begin{equation*}
    h^1\left(\Omega^2_Y\right) = (104 + 9) - 27 = 86 \ .
\end{equation*}
\end{proof}
\end{example}

Actually the numbers of nodes in $\overline{Y}$, of maximally
independent exceptional rational curves in $Y$ and of maximally
independent vanishing cycles in $\widetilde{Y}$ turn out to be
deeply related. This fact characterizes the global change in
topology induced by a conifold transition, as explained in the
following

\begin{theorem}[\cite{Clemens83}, \cite{Reid87}, \cite{Werner-vanGeemen90}, \cite{Tian92},
\cite{Namikawa-Steenbrink95}, \cite{Morrison-Seiberg97},
...]\label{cambio omologico} Let $T(Y,\overline{Y},\widetilde{Y})$
be a conifold transition and let
\begin{itemize}
    \item $N$ be the number of nodes composing $\Sing
    (\overline{Y})$,
    \item $k$ be the maximal number of homologically independent
    exceptional rational curves in $Y$,
    \item $c$ be the maximal number of homologically independent
    vanishing cycles in $\widetilde{Y}$.
\end{itemize}
Then:
\begin{enumerate}
    \item $|\Sing (\overline{Y})|=:N=k+c$;
    \item (Betti numbers) $b_i(Y)=b_i(\overline{Y})=b_i(\widetilde{Y})$ for $i\neq
    2,3,4$, and
    \begin{equation*}
    \begin{array}{ccccc}
      b_2(Y) & = & b_2(\overline{Y})+k & = & b_2(\widetilde{Y})+k \\
      \parallel &  &  &  & \parallel \\
      b_4(Y) & = & b_4(\overline{Y}) & = & b_4(\widetilde{Y})+k \\
       & & & & \\
      b_3(Y) & = & b_3(\overline{Y})-c & = & b_3(\widetilde{Y})-2c \\
    \end{array}
    \end{equation*}
    where vertical equalities are given by Poincar\'{e} Duality;
    \item (Hodge numbers)
    \begin{equation*}
    \begin{array}{ccc}
      h^{2,1}(\widetilde{Y}) & = & h^{2,1}(Y)+c \\
       &  &  \\
      h^{1,1}(\widetilde{Y}) & = & h^{1,1}(Y)-k \\
    \end{array}
    \end{equation*}
\end{enumerate}
\end{theorem}

\begin{remark}
Note that point (2) of the previous statement implies that the
conifold $\overline{Y}$ do not satisfy Poincar\'{e} Duality. The
difference $b_4(\overline{Y}) - b_2(\overline{Y})=k$ is called
\emph{the defect of $\overline{Y}$} \cite{Namikawa-Steenbrink95}.
\end{remark}

\begin{remark}\label{cpx<->ka} Point (3) in theorem \ref{cambio omologico}
has the following geometric interpretation: \emph{a conifold transition
increases complex moduli by the maximal number of homologically
independent vanishing cycles and decreases \ka moduli by the
maximal number of homologically independent exceptional rational
curves}.

\noindent The reader is referred to \ref{cy-moduli} for a deeper
understanding, where the \emph{\cy moduli space's} structure will
be quickly described.
\end{remark}

\begin{proof}[Proof of theorem \ref{cambio omologico}]
Let us denote:
\begin{itemize}
    \item $P:= \Sing(\overline{Y})= \{ p_1, \ldots , p_N \}$, the
    singular locus of $\overline{Y}$;
    \item $E:= \bigcup_{i=1}^{N} E_i$, the exceptional locus of
    $Y$;
    \item $S:= \bigcup_{i=1}^{N} S_i$, the vanishing locus of
    $\widetilde{Y}$.
\end{itemize}
The birational contraction $\phi : Y\rightarrow\overline{Y}$
induces the isomorphism
\begin{equation}\label{resol_iso}
\phi :Y\setminus E \overset{\cong}{\longrightarrow}
\overline{Y}\setminus P
\end{equation}
On the other hand, for any $i=1,\ldots , N$, by  Proposition
\ref{clemens lemma} we can construct  compact tubular
neighborhoods $\widetilde{D}_i$ of the vanishing cycle $S_i$ in
$\widetilde{Y}$ and $\widehat{D}_i$ of the exceptional rational
curve $E_i$ in $Y$ and diffeomorphisms
\begin{equation*}
    \alpha _i:\widehat{D}_i\setminus
    E_i\overset{\cong}{\longrightarrow}\widetilde{D}_i\setminus S_i
\end{equation*}
Since we can clearly assume that $\widetilde{D}_i$'s are all
disjoint neighborhoods and the same for $\widehat{D}_i$'s, the
composed morphisms $\phi\circ\alpha_{i}^{-1}$ give diffeomorphisms
\begin{equation}\label{diffeo-locale}
    \phi\circ\alpha _i^{-1} : \widetilde{D}_i\setminus
    S_i\overset{\cong}{\longrightarrow}\overline{D}_i\setminus \{p_i\}
\end{equation}
where $\overline{D}_i:= \phi\left( \widehat{D}_i\right)$. Set:
\begin{eqnarray*}
  \widetilde{D} &=& \bigcup_{i=0}^{N} \widetilde{D}_i \\
  \overline{D} &=& \bigcup_{i=0}^{N}\overline{D}_i
\end{eqnarray*}
By the Ehreshmann fibration theorem there exists a diffeomorphism
\begin{equation*}
    \widetilde{Y}\setminus
    \widetilde{D}\overset{\cong}{\longrightarrow} \overline{Y}\setminus\overline{D}
\end{equation*}
allowing to extend diffeomorphisms (\ref{diffeo-locale}) to the
following global one
\begin{equation}\label{diffeo-globale}
    \psi : \widetilde{Y}\setminus S \overset{\cong}{\longrightarrow}
    \overline{Y}\setminus P
\end{equation}

\begin{step I}
$\forall i\neq 2,3\ \ b_i(Y)=b_i(\overline{Y})$ and
\begin{equation*}
b_2(Y)= b_2(\overline{Y})+k \Leftrightarrow
b_3(\overline{Y})=b_3(Y)+N-k \ .
\end{equation*}
\end{step I}
Let $T(\widehat{U}_i,\overline{U}_i,\widetilde{U}_i)$ be the local
conifold transition (notation as in section \ref{analisi locale})
induced by $T(Y,\overline{Y},\widetilde{Y})$ around the node
$p_i\in P$ and denote:
\begin{itemize}
    \item $\widehat{U}:=\bigcup_{i=1}^{N} \widehat{U}_i \subset Y$, $Y^*:=
    Y\setminus E$, $\widehat{U}^*:= \widehat{U}\setminus E$;
    \item $\overline{U}:=\bigcup_{i=1}^{N} \overline{U}_i \subset \overline{Y}$, $\overline{Y}^*:=
    \overline{Y}\setminus P$, $\overline{U}^*:= \overline{U}\setminus P$;
\end{itemize}
Then:
\begin{itemize}
    \item $\widehat{U}^* = Y^*\cap \widehat{U}$ and
    $Y=Y^*\cup\widehat{U}$,
    \item $\overline{U}^* = \overline{Y}^*\cap \overline{U}$ and
    $\overline{Y}=\overline{Y}^*\cup\overline{U}$,
\end{itemize}
and we are in a position to apply Mayer--Vietoris machinary to the
couples $(Y^*,\widehat{U})$ and $(\overline{Y}^*,\overline{U})$ to
get the following two long exact sequences in homology
\begin{equation}\label{MVresol}
    \xymatrix@1{\cdots\ar[r] & H_{i}(\widehat{U}^*) \ar[r] &
    H_{i}(Y^*)\oplus H_{i} (\widehat{U}) \ar[r] & H_i (Y) \ar[r] &
    H_{i-1}(\widehat{U}^*)\ar[r] & \cdots}
\end{equation}
\begin{equation}\label{MVsing}
    \xymatrix@1{\cdots\ar[r] & H_{i}(\overline{U}^*) \ar[r] &
    H_{i}(\overline{Y}^*)\oplus H_{i} (\overline{U}) \ar[r] & H_i (\overline{Y}) \ar[r] &
    H_{i-1}(\overline{U}^*)\ar[r] & \cdots}
\end{equation}
By straight line homotopy we have
\begin{equation}\label{slh-resol}
    H_i(\widehat{U}) \cong H_i(E) \cong \left\{ \begin{array}{cc}
      \Z ^N & \text{if}\ i=0,2 \\
      0 & \text{otherwise} \\
    \end{array} \right.
\end{equation}
as a consequence of Proposition \ref{risol-geom} and
\begin{equation}\label{slh-sing}
    H_i(\overline{U}) \cong H_i(P) \cong \left\{ \begin{array}{cc}
      \Z ^N & \text{if}\ i=0 \\
      0 & \text{otherwise} \\
      \end{array} \right.
\end{equation}
as a consequence of Proposition \ref{node-top}. The diffeomorphism
$\phi$, given in (\ref{resol_iso}), induces then the following
isomorphisms in homology
\begin{equation}\label{phi-iso}
    H_i (\widehat{U}^*) \cong H_i(\overline{U}^*) \cong
    \bigoplus_{i=1}^{N}H_i(S^3\times S^2) \cong \left\{ \begin{array}{cc}
      \Z ^N & \text{if}\ i=0,2,3,5 \\
      0 & \text{otherwise} \\
    \end{array} \right.
\end{equation}
and
\begin{equation}\label{phi-iso,bis}
    H_i(Y^*) \cong H_i (\overline{Y}^*) \ .
\end{equation}
Introduce isomorphisms (\ref{slh-resol}), (\ref{slh-sing}),
(\ref{phi-iso}) and (\ref{phi-iso,bis}), as vertical arrows
connecting sequences (\ref{MVresol}) and (\ref{MVsing}). The
Steenrod 5--lemma gives then
\begin{equation}\label{ts 1}
    \forall i\neq 2,3 \quad  b_i (Y) = b_i (\overline{Y})\ .
\end{equation}
Moreover gluing the two sequences by identifying the isomorphic
poles, they reduce to the following diagram
\begin{eqnarray}\label{MVdiag1}
  && \xymatrix@1{0\ar[r]& H_4(Y^*)\ar[r] & H_4(Y)\ar[r] & \cdots } \\
  \nonumber &&\xymatrix{&H_3(Y)\ar[rd]&&H_2(Y^*)\oplus\Z ^N\ar[r]&H_2(Y)\ar[rd]& \\
              \cdots H_3(Y^*)\ar[ru]\ar[rd]&&\Z ^N\ar[ru]\ar[rd]&&&0 \\
              &H_3(\overline{Y})\ar[ru]&&H_2(\overline{Y}^*)\ar[r]&H_2(\overline{Y})\ar[ru]& }
\end{eqnarray}
Then we get the following relations on Betti numbers
\begin{eqnarray*}
  b_4(Y^*)-b_4(Y)+b_3(Y^*)-b_3(Y)+N-(b_2(Y^*)+N) + b_2(Y) &=& 0 \\
  b_4(Y^*)-b_4(Y)+b_3(Y^*)-b_3(\overline{Y})+N-b_2(\overline{Y}^*) + b_2(\overline{Y}) &=& 0
\end{eqnarray*}
and their difference gives
\begin{equation*}
    b_2(Y) - b_2(\overline{Y}) = b_3(Y)  - b_3(\overline{Y})+ N\ .
\end{equation*}

\begin{step II}
$\forall i\neq 3,4\ \ b_i(\widetilde{Y})=b_i(\overline{Y})$ and
\begin{equation*}
b_3(\widetilde{Y})= b_3(\overline{Y})+c \Leftrightarrow
b_4(\overline{Y})=b_4(\widetilde{Y})+N-c \ .
\end{equation*}
\end{step II}

Let $T(\widehat{U}_i,\overline{U}_i,\widetilde{U}_i)$ be the local
conifold induced near the node $p_i \in P$, as before. Let us
denote
\begin{itemize}
    \item $\widetilde{U}:=\bigcup_{i=1}^{N}
    \widetilde{U}_i\subset\widetilde{Y}$, $\widetilde{Y}^*:=\widetilde{Y}\setminus
    S$, $\widetilde{U}^*:=\widetilde{U}\setminus S$.
\end{itemize}
Then
\begin{itemize}
    \item $\widetilde{U}^* = \widetilde{Y}^* \cap\widetilde{U}$ and
    $\widetilde{Y} = \widetilde{Y}^* \cup\widetilde{U}$
\end{itemize}
and Maeyer--Vietoris sequence for the couple
$(\widetilde{Y}^*,\widetilde{U})$ gives
\begin{equation}\label{MVsmooth}
    \xymatrix@1{\cdots \ar[r]&H_i(\widetilde{U}^*)\ar[r]&H_i(\widetilde{Y}^*)\oplus H_i(\widetilde{U})\ar[r]&
    H_i(\widetilde{Y})\ar[r]&H_{i-1}(\widetilde{U}^*)\ar[r]&\cdots  }
\end{equation}
Proposition \ref{smooth-geom} and straight line homotopy give
\begin{equation}\label{slh-smooth}
    H_i (\widetilde{U}) \cong H_i(S) \cong
    \bigoplus_{i=1}^{N}H_i(S^3) \cong \left\{ \begin{array}{cc}
      \Z ^N & \text{if}\ i=0,3, \\
      0 & \text{otherwise} \\
    \end{array} \right.
\end{equation}
Moreover the diffeomorphism $\psi$ given in (\ref{diffeo-globale})
induces the following isomorphisms in homology
\begin{equation}\label{psi-iso}
    H_i (\widetilde{U}^*) \cong H_i(\overline{U}^*) \cong
    \bigoplus_{i=1}^{N}H_i(S^3\times S^2) \cong \left\{ \begin{array}{cc}
      \Z ^N & \text{if}\ i=0,2,3,5 \\
      0 & \text{otherwise} \\
    \end{array} \right.
\end{equation}
and
\begin{equation}\label{psi-iso,bis}
    H_i(\widetilde{Y}^*) \cong H_i (\overline{Y}^*) \ .
\end{equation}
As before, apply the Steenrod 5--lemma to conclude that
\begin{equation}\label{ts 3}
    \forall i\neq 3,4 \quad  b_i (\widetilde{Y}) = b_i (\overline{Y})\ .
\end{equation}
and glue sequences (\ref{MVsing}) and (\ref{MVsmooth}) to get the
following diagram
\begin{eqnarray}\label{MVdiag2}
  &&\xymatrix{&&H_4(\widetilde{Y})\ar[rd]&&H_3(\widetilde{Y}^*)\oplus\Z ^N\ar[r]&H_3(\widetilde{Y})\cdots \\
              0\ar[r] &H_4(\widetilde{Y}^*)\ar[ru]\ar[rd]&&\Z ^N\ar[ru]\ar[rd]&& \\
              &&H_4(\overline{Y})\ar[ru]&&H_3(\overline{Y}^*)\ar[r]&H_3(\overline{Y})\cdots}\\
  \nonumber &&\xymatrix{\cdots\ar[rd]&&&&\\
                        &\Z ^N\ar[r]& H_2(\widetilde{Y}^*)\ar[r]& H_2(\widetilde{Y})\ar[r]&0\\
                        \cdots\ar[ru]&&&& }
\end{eqnarray}
Then we get the following relations on Betti numbers
\begin{eqnarray*}
  b_4(\widetilde{Y}^*)-b_4(\widetilde{Y})+N-(b_3(\widetilde{Y}^*)+N)+b_3(\widetilde{Y})-N+b_2(\widetilde{Y}^*)-
  b_2(\widetilde{Y}) &=& 0 \\
  b_4(\widetilde{Y}^*)-b_4(\overline{Y})+N-b_3(\widetilde{Y}^*)+b_3(\overline{Y})-N+b_2(\widetilde{Y}^*)-
  b_2(\widetilde{Y}) &=& 0
\end{eqnarray*}
and their difference gives
\begin{equation*}
    b_3(\widetilde{Y}) - b_3(\overline{Y}) = b_4(\widetilde{Y}) - b_4(\overline{Y})+ N\ .
\end{equation*}

\begin{step III}
Let $k$ and $c$ be the same parameters defined in Steps I and II
respectively. Then
\begin{equation*}
    |\Sing (\overline{Y})|=:N=k+c:=
    b_2(Y)-b_2(\overline{Y})+b_3(\widetilde{Y})-b_3(\overline{Y})\
    .
\end{equation*}
\end{step III}

By Poincar\'e duality
\begin{eqnarray*}
  b_2(Y) &=& b_4(Y) \\
  b_4(\widetilde{Y}) &=& b_2(\widetilde{Y})
\end{eqnarray*}
Recall then Steps I and II to get
\begin{eqnarray*}
b_2(Y)&=&b_4(Y)=b_4(\overline{Y})=b_4(\widetilde{Y})+N-c\\
&=&b_2(\widetilde{Y})+N-c=b_2(\overline{Y})+N-c=b_2(Y)-k+N-c\\
\end{eqnarray*}
Hence $N-k-c=0$.

\begin{step IV}
$k$ is the maximal number of homologically independent exceptional
rational curves in $Y$ while $c$ is the maximal number of
homologically independent vanishing cycles in $\widetilde{Y}$.
\end{step IV}

Recall the diffeomorphisms $\phi$ and $\psi$, defined in
(\ref{resol_iso}) and (\ref{diffeo-globale}),  and consider the
composition
\begin{equation}\label{resol-smooth_iso}
    \psi^{-1}\circ\phi: Y\setminus E
    \overset{\cong}{\longrightarrow} \widetilde{Y}\setminus S
\end{equation}
Lefschetz duality ensures that
\begin{eqnarray*}
  H^{6-i}(Y\setminus E) &\cong& H_i (Y,E) \\
  H^{6-i}(\widetilde{Y}\setminus S) &=& H_i (\widetilde{Y},S)
\end{eqnarray*}
Then (\ref{resol-smooth_iso}) gives
\begin{equation}\label{rel-hom_iso}
H_i (Y,E) \cong H_i (\widetilde{Y},S)
\end{equation}
Consider the long exact relative homology sequences of the couples
$(Y,E)$ and $(\widetilde{Y},S)$ and the vertical isomorphisms
given by (\ref{rel-hom_iso}):
\begin{equation}\label{rel-hom_diagram}
    \xymatrix{\cdots H_{i+1}(Y,E)\ar[r]\ar[d]^{\cong}& H_i(E)\ar[r]& H_i(Y)\ar[r]& H_i(Y,E)\cdots\ar[d]^{\cong}\\
              \cdots H_{i+1}(\widetilde{Y},S)\ar[r]& H_i(S)\ar[r]& H_i(\widetilde{Y})\ar[r]& H_i(\widetilde{Y},S)\cdots}
\end{equation}
By identifying the isomorphic poles and recalling
(\ref{slh-resol}) and (\ref{slh-smooth}) the previous long exact
sequences reduce to the following diagram:
\begin{equation}\label{k-c-diagram}
    \xymatrix{&&&&&0\ar[d]&\\
              &&&&&H_3(Y)\ar[d]&\\
              0\ar[r]&H_4(\widetilde{Y})\ar[r]&H_4(Y)\ar[r]&H_3(S)\ar[r]^{\gamma}\ar@{}[d]|{\parallel}&H_3(\widetilde{Y})\ar[r]&
              H_3(\widetilde{Y},S)\ar[r]\ar[d]&0\\
              &&&\Z ^N &&H_2(E)\ar@{}[r]|{=}\ar[d]^{\kappa}&\Z ^N\\
              &&&&&H_2(Y)\ar[d]&\\
              &&&&&H_2(\widetilde{Y})\ar[d]&\\
              &&&&&0&}
\end{equation}
Set
\begin{equation*}
    I:=\im [\kappa :\Z ^N = H_2(E)\longrightarrow H_2(Y)]
\end{equation*}
Then $k:= \rk (I)$ is the number of linear independent classes of
exceptional curves in $H_2(Y)$. Since
\begin{equation*}
    \xymatrix@1{0\ar[r]&I\ar[r]&H_2(Y)\ar[r]&H_2(\widetilde{Y})\ar[r]&0}
\end{equation*}
is a short exact sequence, it follows that
\begin{equation*}
    b_2(Y)=b_2(\widetilde{Y})+k
\end{equation*}
On the other hand set
\begin{equation*}
    K:=\ker [\gamma :\Z ^N\cong H_3(S)\longrightarrow H_3(\widetilde{Y})]
\end{equation*}
Then $N-c:=\rk (K)$ is the number of linear independent relations
on the classes of vanishing cycles in $H_3(\widetilde{Y})$. Since
\begin{equation*}
    \xymatrix@1{0\ar[r]&H_4(\widetilde{Y})\ar[r]&H_4(Y)\ar[r]&K\ar[r]&0}
\end{equation*}
is a short exact sequence, it follows that
\begin{equation*}
    b_4(Y)=b_4(\widetilde{Y})+N-c
\end{equation*}
\end{proof}

\subsection{What about more general geometric
transitions?}\label{cambio omologico generalizzato}

The local and global topology and geometry of a general geometric
transition
\begin{equation*}
    \xymatrix@1{Y\ar@/_1pc/ @{.>}[rr]_T\ar[r]^{\phi}&
                \overline{Y}\ar@{<~>}[r]&\widetilde{Y}}\ .
\end{equation*}
can actually be very intricate, depending on the nature of
$\Sing(\overline{Y})$ and on the geometry of the exceptional locus
of $\phi$. For this reason no general results similar to
Proposition \ref{clemens lemma} and theorem \ref{cambio omologico}
are known. Anyway, under some (strong) condition on
$\Sing(\overline{Y})$, somewhat can be said.

\noindent First of all let us assume that
\emph{$\Sing(\overline{Y})=\{ p_1,\ldots,p_r\}$ is composed only
by isolated hypersurface singularity}.

\noindent In this case, given a 1--parameter flat smoothing
$\mathcal{Y}\rightarrow\Delta^1$ of the singular point $p_i$, the
local topology of $\mathcal{Y}$ near $p_i$ is explained by the
Milnor's analysis \cite{Milnor68}. Call $B$ the union of all of
the Milnor's fibres $B_{p_i}$, which have the homology type of a
bouquet of 3--spheres. Interpolate the relative homology long
exact sequences of $(\overline{Y},\Sing(\overline{Y}))$ and
$(\widetilde{Y},B)$, like in step IV of the proof of theorem
\ref{cambio omologico}, to get the first part of the following

\begin{theorem}[\cite{Namikawa-Steenbrink95}, theorem
(3.2)]\label{NS-thm} Let $\overline{Y}$ be a normal projective
3--fold with only isolated hypersurface singularities, admitting a
smoothing $\widetilde{Y}$. For any $p\in\Sing(\overline{Y})$ call
$m(p):=h_3(B_p)$ the \emph{Milnor number} of $p$. Then the defect
of $\overline{Y}$ is related to Milnor numbers as follows
\begin{equation}\label{cambio omologico II}
    k:= b_4 \left(\overline{Y}\right) - b_2 \left(\overline{Y}\right)
    = b_3 \left(\overline{Y}\right) + \sum_{p\in\Sing(\overline{Y})}
    m(p) - b_3 \left(\widetilde{Y}\right) \ .
\end{equation}
Moreover if all of the singularities of $\overline{Y}$ are
rational then
\begin{equation*}
    W\left(\overline{Y}\right)/C\left(\overline{Y}\right):=
    \left\langle\text{Weil divisors of }\overline{Y}\right\rangle_{\Z} /
    \left\langle\text{Cartan divisors of }\overline{Y}\right\rangle_{\Z}
\end{equation*}
is a finitely generated abelian group. In particular if
$h^2(\mathcal{O}_{\overline{Y}})=0$ then
\begin{equation*}
    k=\rk \left( W\left(\overline{Y}\right)/C\left(\overline{Y}\right)\right)
\end{equation*}
giving a further interpretation of the defect of $\overline{Y}$.
\end{theorem}

Since the Milnor fibre of a node $p$ has the homology type of a
single 3--sphere, $m(p)=1$ and (\ref{cambio omologico II}) gives
(1) and the right part of formulas (2) in theorem \ref{cambio
omologico}.

\noindent The last part of the previous statement is proved by
employing results of A.~Dimca \cite{Dimca90} and J.~H.~M.~Steenbrink
\cite{Steenbrink85}.

\noindent Moreover theorem \ref{NS-thm}, joint with results of
M.~Reid \cite{Reid83}, allows to generalize theorem \ref{cambio
omologico} to the case of a geometric transition whose birational
contraction is a small one, as follows.

\begin{theorem}[\cite{Namikawa-Steenbrink95}, Example (3.8)]\label{NS-ex}
Let $T(Y,\overline{Y},\widetilde{Y})$ be a geometric transition
whose birational contraction $\phi:Y\rightarrow\overline{Y}$ is a
composition of type I primitive contractions. Then
$\Sing{\overline{Y}}=\{ p_1,\ldots,p_r\}$ where $p_i$ is an
isolated, rational singularity. Let $C_i:=\phi ^{-1}(p_i)$ and
$n_i$ be the number of irreducible components of $C_i$. If $k$ is
the rank of the free abelian group generated in $H_2(Y)$ by the
homology classes of $C_1,\ldots,C_r$, then
\begin{eqnarray*}
    b_2(\widetilde{Y})&=&b_2(Y)-k\\
    b_3(\widetilde{Y})&=&b_3(Y)+\sum _{i=1}^r n_i + \sum _{i=1}^r m(p_i) -2k
\end{eqnarray*}
\end{theorem}

As far as I know, dropping assumptions on $\Sing(\overline{Y})$
leads to no more than interesting conjectures and examples. The
interested reader is referred to \cite{Morrison-Seiberg97},
section 3 and appendix A, for some geometric and physical
interpretation of parameters $N,k,c$ for more general transitions,
and to \cite{KMP96} for a computation of these parameters in
examples of transitions whose $\overline{Y}$ admits non--isolated
singularities (see also \ref{transizione reverse} in the
following).

\section{Classification of geometric transitions}

By definition, a general geometric (not necessarily conifold)
transition $T(Y,\overline{Y},\widetilde{Y})$ is always associated
with a \emph{birational contraction} of a \cy threefold $Y$ to a
normal variety $\overline{Y}$. Then the ingredients of a
\emph{classification} are the following:
\begin{enumerate}
    \item to classify the birational contractions
    $\phi:Y\rightarrow\overline{Y}$ which may occur,
    \item among them, to select those admitting a smoothable
    target $\overline{Y}$.
\end{enumerate}

\noindent Let us start with the first point of our program.

\subsection{A little bit of Mori theory for \cy threefolds}

Let $Y$ be a \cy threefold and consider the \emph{Picard group}
\begin{eqnarray*}
  \Pic (Y) &:=& \left\langle \text{Invertible Sheaves}\right\rangle_{\Z} /\text{isomorphism ($\cong$)} \\
   &\cong& \left\langle \text{Divisors}\right\rangle_{\Z} /\text{linear equivalence ($\equiv$)}
\end{eqnarray*}

\begin{remark}
There is a canonical isomorphism
\begin{equation}\label{Pic-cy}
    \Pic (Y)\cong H^2(Y,\Z)
\end{equation}
In fact, since $Y$ is smooth, it is well known that $\Pic (Y)\cong
H^1(Y,\mathcal{O}_Y^*)$. The long exact cohomology sequence
associated with the \emph{exponential sequence}
\begin{equation*}
    \xymatrix@1{0\ar[r]& \Z \ar[r]& \mathcal{O}\ar[r]^{\text{exp}}& \mathcal{O}^*\ar[r]& 1}
\end{equation*}
gives the claim as a consequence of the \cy condition
$h^1(\mathcal{O}_Y)=h^2(\mathcal{O}_Y)=0$.
\end{remark}

The \emph{Kleiman space} is the following real vector space
\begin{equation}\label{Kleiman-vs}
    H^2(Y,\R)\cong H^2(Y,\Z)\otimes_{\Z}\R \cong \R ^{\rho}
\end{equation}
whose dimension is clearly $\rho = \rk (\Pic (Y))$, called the
\emph{Picard number of $Y$}.

\begin{definition}
A divisor $D$ of $Y$ is called \emph{nef} (\emph{numerically
effective}) if any curve $C$ in $Y$ intersects $D$ non--negatively
i.e.
\begin{equation*}
    (D\cdot C)\geq 0
\end{equation*}
\end{definition}

\begin{definition}[The closed \ka cone]\label{cono di kahler}
The \emph{closed \ka cone} $\overline{\mathcal{K}}(Y)$ of $Y$ is
the cone generated in the Kleiman space $H^2(Y,\R)$ by the classes
of nef divisors.
\end{definition}

\begin{definition}[The closed Mori cone]
The dual construction with respect to the perfect pairing
\begin{equation*}
    (\ \cdot\ ):H^2(Y,\R)\otimes H_2(Y,\R)\longrightarrow \R
\end{equation*}
induced by the intersection product, gives rise to the
\emph{closed Mori cone} $\overline{NE}(Y)$.
\end{definition}

\begin{theorem}[Kleiman Ampleness Criterion
\cite{Kleiman66}]\label{KAC} A divisor $D$ of $Y$ (not necessarily
neither \cy nor 3--dimensional) is ample if and only if
\begin{equation*}
    \forall Z\in \overline{NE}(Y)\setminus\{ 0\} \quad (D\cdot
    Z)>0
\end{equation*}
\end{theorem}

\begin{corollary}\label{kahler-interpretazione}
Let $Y$ be \cy variety. The interior $\mathcal{K}(Y)$ of
$\overline{\mathcal{K}}(Y)$ is the cone generated by the \ka
classes in the Kleiman space $H^2(Y,\R)$.
\end{corollary}

\begin{proof}
The criterion \ref{KAC} ensures that $\mathcal{K}(Y)$ is the cone
generated by the classes of ample divisors in $H^2(Y,\R)$. A
divisor is ample if and only if its fundamental form is positive,
then $D$ is ample if and only if $[D]\in H^2(Y,\R)$ is the class
of a \ka form, since the \cy condition ensures that
$H^2(Y,\C)\cong H^{1,1}(Y)$.
\end{proof}

\begin{theorem}[of the Mori cone \cite{Mori82}]
The \emph{negative part} of $\overline{NE}(Y)$ ($Y$ not
necessarily neither \cy nor 3--dimensional) is rational and
polyhedral i.e. there exists a collection $\{ C_i\}_{i\in I}$ of
rational curves in $Y$ such that
\begin{equation*}
\overline{NE}(Y)_{-}:=\overline{NE}(Y)\cap \{ Z\in
\overline{NE}(Y)|(K_Y\cdot Z)<0\}=\sum_{i\in I} \R _{\geq 0}[C_i]
\ .
\end{equation*}
\end{theorem}

\begin{theorem}[of the \ka cone \cite{Wilson89}, \cite{Wilson92}]
Let $Y$ be a \cy threefold and consider the cubic cone in
$H^2(Y,\R)$ given by the cup--product
\begin{equation}\label{W*}
    W^*:=\{ [D]\in H^2(Y,\R)|D^3=0\}
\end{equation}
(it is the cone projecting a cubic hypersurface $W\subset \P
(H^2(Y,\R))=\P ^{\rho-1}_{\R}$). Then
\begin{equation}\label{W*_in_bd}
    W^* \cap \overline{\mathcal{K}}(Y) \subset\partial \overline{\mathcal{K}}(Y)
\end{equation}
and $\overline{\mathcal{K}}(Y)$ is locally polyhedral away from
$W^*$. In particular $\partial \overline{\mathcal{K}}(Y)\setminus
W^*$ is composed by codimension 1 faces and their intersections.
\end{theorem}

\begin{remark}
(\ref{W*_in_bd}) is an immediate consequence of the definition of
$\overline{\mathcal{K}}(Y)$. In fact if there exists $[D]\in
W^*\cap \mathcal{K}(Y)$ then $D$ should be ample, implying that
$D^3>0$ and contradicting (\ref{W*}).
\end{remark}

\begin{remark}
By Corollary \ref{kahler-interpretazione},
$\partial\overline{\mathcal{K}}(Y)\cap W^*$ parameterizes all the
possible \emph{degenerations} of a \ka metric on $Y$ (see
\cite{Morrison-Seiberg97}, section 3).
\end{remark}

\subsection{Primitive contractions and primitive transitions}

\begin{definition}
Let $\phi:Y\rightarrow\overline{Y}$ be a birational contraction of
a \cy variety to a normal one. $\phi$ is called \emph{primitive}
(or alternatively \emph{extremal}, as explained in remark
\ref{1:1}.(1)) if it cannot be factored into birational morphisms
of normal varieties. Any associated transition
$T(Y,\overline{Y},\widetilde{Y})$ is called a \emph{primitive} (or
\emph{extremal}) {transition}.
\end{definition}

\begin{proposition}[Contractions by the Mori--\ka cones point of
view]\label{corrispondenza}
There is a correspondence
\begin{equation*}
    \left\{ \begin{array}{c}
      \phi : Y\rightarrow\overline{Y}\ \text{contraction} \\
      \text{from \cy to normal} \\
    \end{array}\right\}\leftrightarrow
    \left( \partial\overline{\mathcal{K}}(Y)\setminus W^*\right)_{\Q}
    \leftrightarrow \left( \partial\overline{NE}(Y)\cap \overline{NE} (Y)_{-}\right)_{\Q}
\end{equation*}
where $(\ )_{\Q}$ means ``rational points of''. In particular
\begin{eqnarray*}
  \phi\ \text{is primitive} &\Leftrightarrow& \begin{array}{c}
    \text{it corresponds to a class $[D]$ in the interior} \\
    \text{of a codimension 1 face of $\overline{\mathcal{K}}(Y)$} \\
  \end{array} \\
   &\Leftrightarrow& \begin{array}{c}
     \text{it corresponds to a class} \\
     \text{generating an extremal ray of $\overline{NE}(Y)$} \\
   \end{array}
\end{eqnarray*}
\end{proposition}

\begin{proof}[Sketch of proof]
Let $H$ be a hyperplane section of $\overline{Y}$. Since
$\overline{Y}$ is normal we can assume
\begin{equation}\label{normalita'}
    H\cap \Sing (\overline{Y})=\emptyset\ .
\end{equation}
Look at the pull--back $\phi^* H$. The Kleinman Criterion
\ref{KAC} ensures that
\begin{equation*}
    \forall Z \in \overline{NE}(Y) \quad (\phi ^* H \cdot Z)\geq 0
\end{equation*}
In particular, if $E$ is the exceptional locus of $\phi$, the
\emph{projection formula} and (\ref{normalita'}) give
\begin{equation*}
(\phi ^* H \cdot Z)=0 \Leftrightarrow \text{$Z$ is the class of a
curve $C\subset E$}\ .
\end{equation*}
Then $(\phi ^* H\cdot\quad)$ defines a hyperplane in $H_2(Y,\R)$
cutting $\overline{NE}(Y)$ along an extremal face. By duality
$[\phi ^* H]$ generates a ray living in a codimension 1 face of
the polyhedral part of the \ka cone i.e.
\begin{equation*}
    \R_{\geq 0}[\phi ^* H]\subset
    \partial\overline{\mathcal{K}}(Y)\setminus W^*\ .
\end{equation*}

\noindent Notice that the contraction $\phi$ can be factored into
birational morphisms if there exists a curve $C$ in $E$ and
$Z_1,Z_2\in \overline{NE}(Y)$ such that
\begin{equation}
    \R_{\geq 0}Z_1\neq \R_{\geq 0}Z_2\quad \text{and}\quad [C]=Z_1+Z_2\ .
\end{equation}
Hence
\begin{eqnarray*}
  \phi \ \text{is primitive} &\Leftrightarrow& \forall C\subset E\quad \R _{\geq 0}[C]\quad
  \text{is the same extremal ray of}\ \overline{NE}(Y) \\
   &\Leftrightarrow& \begin{array}{c}
     \R_{\geq 0}[\phi ^* H] \quad \text{is not on the intersection} \\
     \text{of two codimension 1 faces of}\ \overline{\mathcal{K}}(Y) \\
    \end{array}\\
   &\Leftrightarrow& \phi ^* H \ \text{is an interior point of an extremal cod. 1 face}\\
\end{eqnarray*}
\end{proof}

\begin{corollary}
Let $T(Y,\overline{Y},\widetilde{Y})$ be a geometric transition
and $\phi :Y\rightarrow\overline{Y}$ the associated birational
contraction. Then $\phi$ can always be factored into a composite
of a \emph{finite} number of primitive contractions.
\end{corollary}

\begin{remark}
The finiteness of the factorization process follows from the fact
that any primitive contraction reduces by 1 the Picard number.
\end{remark}

\begin{remark}\label{1:1}
The correspondence given in Proposition \ref{corrispondenza} from
contraction morphisms and rational points of the boundary of the
\ka cone \emph{is not a 1:1 correspondence}. Actually all the
rational classes living in the interior of the same codimension 1
face of the $\overline{\mathcal{K}}(Y)$ correspond to the same
primitive birational contraction.

\noindent It is then possible to conclude that:
\begin{enumerate}
    \item \emph{there is a 1:1 correspondence between primitive contractions and either codimension 1 faces
    of the \ka cone $\overline{\mathcal{K}}(Y)$ or extremal rays of the Mori cone
    $\overline{NE}(X)$} (\cite{Wilson92}, fact 1); for this reason primitive contractions (transitions) are also
    called \emph{extremal} contractions (transitions) \cite{Morrison99};
    \item \emph{there is a 1:1 correspondence between codimension $r$ faces of the
    \ka cone $\overline{\mathcal{K}}(Y)$ and birational contractions from a \cy 3--fold
    to a normal variety composed by $r$ primitive contractions}.
\end{enumerate}
\end{remark}

\begin{theorem}[Classification of primitive contraction
\cite{Wilson92}]\label{classificazione} Let $\phi :
Y\rightarrow\overline{Y}$ be a primitive contraction from a \cy
threefold to a normal variety. Then one of the following is true:
\begin{description}
    \item[type I] $\phi$ is \emph{small} and the exceptional locus
    $E$ is composed of finitely many rational curves;
    \item[type II] $\phi$ contracts a divisor down to a point; in
    this case $E$ is irreducible and in particular it is a
    (generalized) \emph{del Pezzo surface} (see \cite{Reid80})
    \item[type III] $\phi$ contracts a divisor down to a curve
    $C$; in this case $E$ is still irreducible and it is a conic
    bundle over a smooth curve $C$.
\end{description}
\end{theorem}

\begin{definition}[Classification of primitive transitions]
A transition $T(Y,\overline{Y},\widetilde{Y})$ is called \emph{of
type I, II or III} if it is \emph{primitive} and if the associated
birational contraction $\phi : Y\rightarrow\overline{Y}$ is of
type I, II or III, respectively.
\end{definition}

\subsection{Smoothing the target space $\overline{Y}$}

Let us now consider the second point of the classification program
given at the beginning of the present section.

\noindent Let $\phi : Y\rightarrow\overline{Y}$ be a birational
contraction of a \cy 3--fold to a normal one. The problem is to
select all those contractions admitting a smoothable target space
$\overline{Y}$.

\noindent To answer need to analyze the singularities of
$\overline{Y}$ and actually the geometry of the exceptional locus
of $\phi$. Since this is a very hard (and almost completely open)
problem for a general birational contraction $\phi$ let us at
first restrict to consider the case of primitive contractions, as
classified by theorem \ref{classificazione}.

\subsubsection{Transitions of type I}\label{tipo I}

$\phi$ is the contraction of $E_1,\ldots ,E_N$ with $E_i\cong
\P^1$. Then:
\begin{enumerate}
    \item $\overline{Y}$ has $N$ isolated singularities $p_i =
    \phi (E_i)$.
    \item Reid proved that isolated singularities of this kind are actually
    \emph{compound Du Val (cDV) singularities} (see
    \cite{Reid83}, Corollary (1.12)) i.e. they admit local equation of the following
    type
    \begin{equation}\label{cdV}
    f(x,y,z)+tg(x,y,z,t)=0 \quad \text{in $\C^4$}
    \end{equation}
    where $f(x,y,z)=0$ is the local equation in $\C ^3$ of a \emph{rational surface
    singularity} (also known as \emph{Du Val singularity}, see
    \cite{Reid80}, \cite{BPvdV84}). The equation (\ref{cdV}) actually means that
    our 3--dimensional singularity reduces to a rational surface
    singularity on a suitable section.
    \item If $Y$ is \emph{general} (in its complex moduli space)
    such a singular point can be reduced to be an ordinary
    double point (a \emph{node}) i.e.
    \begin{equation*}
    f(x,y,z)=x^2+y^2+z^2\quad\text{and}\quad g(x,y,z,t)=t
    \end{equation*}
    This fact follows from the following:

    \begin{theorem}[\cite{Namikawa94} Theorem B,
    \cite{Namikawa-Steenbrink95} Theorem 2.4,
    \cite{Namikawa02} Theorem 2.7, \cite{Gross97a} Theorem 3.8]\label{ng-thm}
    Let $Y$ be a \cy 3--fold and suppose $\phi:
    Y\rightarrow\overline{Y}$ is a birational contraction morphism
    such that $\overline{Y}$ has isolated, canonical, complete
    intersection singularities. Then there is a deformation of
    $\overline{Y}$to a variety with at worst ordinary double
    points.
    \end{theorem}

    \begin{remark}
    By the previous point (2) we simply may assume $Y$ to admit
    isolated cDV singular points which are, in particular,
    hypersurfaces singularities. Then for the present purpose it suffices Theorem 2.4 of
    \cite{Namikawa-Steenbrink95} to conclude.

    \noindent Anyway we preferred to state Theorem \ref{ng-thm} in
    the improved form given by M.~Gross (\cite{Gross97a} Theorem
    3.8) for further applications in the case of more general
    transitions.
    \end{remark}

    \begin{remark}
    In \cite{Gross97a}, Corollary 3.10, $Y$ may be also assumed
    to be \emph{$\Q$--factorial} (i.e. $\rk (W(Y)/C(Y)) =0$, see Theorem
    \ref{NS-thm})
    \emph{with terminal singularities}. In fact, by results of Y.~Namikawa and J.~Steenbrink
    \cite{Namikawa-Steenbrink95}, \cite{Namikawa94}, in this case there are small
    deformations $\mathcal{Y}\rightarrow\Delta$ and $\overline{\mathcal{Y}}\rightarrow\Delta$ of
    $Y$ and $\overline{Y}$, respectively, such that the morphism $\phi: Y\rightarrow\overline{Y}$ can be deformed
    to a morphism $\varphi :\mathcal{Y}\rightarrow\overline{\mathcal{Y}}$. In particular, for $t\neq
    0$, $\mathcal{Y}_t$ is smooth and $\overline{\mathcal{Y}}_t$
    still has isolated complete intersection singularities but
    admits a crepant resolution $\varphi
    _t:\mathcal{Y}_t\rightarrow\overline{\mathcal{Y}}_t$. Then one
    applies Theorem \ref{ng-thm} to $\overline{\mathcal{Y}}_t$.
    \end{remark}

    \item The last step is the following result essentially due to
    R.~Friedman:

    \begin{theorem}[\cite{Friedman86}, \cite{Friedman91}, \cite{Gross97a} Theorem
    5.1]\label{friedman-thm}
    If $\phi : Y\rightarrow\overline{Y}$ is of type I and
    $\overline{Y}$ has at most ordinary double points then
    $\overline{Y}$ admits a \cy smoothing $\widetilde{Y}$ except
    for the case $N=1$ (\emph{which is:} if $\phi$ contracts a single $\P ^1$
    to a node then $\overline{Y}$ is rigid).
    \end{theorem}

    \begin{proof}[Sketch of proof]
    The key fact in proving the previous theorem is that the
    exceptional curves $E_1,\ldots ,E_N$ of $\phi$ must be \emph{homologically
    dependent} in $H_2(Y,\Z)$, since $\phi$ is a primitive
    contraction i.e. it is the contraction of a unique extremal
    ray $\R _{\geq 0}[E_i]\subset \overline{NE}(Y)$. Then there is
    a non--trivial linear dependence relation on
    $[E_1],\ldots,[E_N]$, except for $N=1$. Results of R.~Friedman, Y.~Namikawa and G.~Tian
    conclude the proof (see \cite{Friedman86} remark 4.5,
    \cite{Friedman91} Proposition 8.7, \cite{Tian92} Theorem 0.1, \cite{Namikawa02} Theorem 2.5).
    \end{proof}

\end{enumerate}

\begin{conclusion}
If $T(Y,\overline{Y},\widetilde{Y})$ is a type I transition then
the exceptional locus $E$ is composed by $N\geq 2$ rational
curves. Moreover if $Y$ is general then $T$ is a conifold
transition contracting $N\geq 2$ rational curves down to nodes.

\noindent In particular $E$ can never be isomorphic to a single
$\P^1$.
\end{conclusion}

\subsubsection{Transitions of type II}

$\phi$ is the contraction of an irreducible divisor $E$ which is a
generalized del Pezzo surface. Then:
\begin{enumerate}
    \item $\overline{Y}$ has one singular point $p=\phi(E)$, which
    is a \emph{canonical singularity} \cite{Reid80},
    \cite{Reid87a}. In particular $\phi$ \emph{is the blowing up of $\overline{Y}$
    at $p$ and the exceptional surface $E$ is either a normal, rational, del Pezzo surface
    of degree $k\leq 9$ or a non--normal del Pezzo surface as classified in
    \cite{Reid94}}.
    \item $k=\deg E$ is the \emph{Reid's invariant of the singularity
    $p=\phi(E)$}. In particular we get that (see \cite{Reid80}, Proposition 2.9 and Corollary 2.10):
    \begin{description}
        \item[$k\leq 2$] then $p$ is a hypersurface singularity whose
        local equation is known,
        \item[$k\geq 3$] then $p$ is a singularity of multiplicity
        $k$ and minimal embedding dimension $\dim \left( m_p/m_p^2\right) = k+1$.
    \end{description}
    In particular, for $k\leq 4$, $p$ \emph{is a complete intersection
    singularity} and, on the contrary, for $k\geq 5$, $p$ \emph{is never a complete intersection
    singularity}.

    \noindent We can then apply Theorem \ref{ng-thm} to conclude
    that \emph{there exists a smoothing $\widetilde{Y}$ of $\overline{Y}$ when $E$ is normal and
    $\deg E \leq 4$}, since $p$ can never be a node.
    \item \emph{When $E$ is normal and $k\geq 5$ then $E$ is smooth and $p$ is
    analytically isomorphic to the vertex of a cone over $E$} (\cite{Gross97a}, Proposition
    5.4). The deformation theory of a cone over a smooth del Pezzo
    surface of degree $5\leq k\leq 9$ is known and precisely:
    \begin{description}
        \item[$k=5$] then $p$ is a codimension 3 singularity and
        there exists a smoothing $\widetilde{Y}$ of $\overline{Y}$
        since locally $\overline{Y}$ is a Pfaffian subscheme
        \cite{Kleppe-Laksov80},
        \item[$6\leq k\leq 9$] then the considered cones are toric
        varieties and by \cite{Altmann97} we get:
        \begin{description}
            \item[$k=6$] then \emph{there are two distinct smoothings}
            $\widetilde{Y}$ given either by the generic
            hyperplane section of a cone over
            $\P^1\times\P^1\times\P^1\subset\P^7$ or by two generic
            hyperplane sections of a cone over
            $\P^2\times\P^2\subset\P^8$;
            \item[$k=7$] then \emph{there is a smoothing} $\widetilde{Y}$
            given by the generic hyperplane section of a cone over
            $\P^3$ blown up at a point, suitably embedded in
            $\P^8$;
            \item[$k=8$] then either $E\cong \P^1\times\P^1$ and \emph{there exists a
            smoothing} $\widetilde{Y}$ given by the generic
            hyperplane section of a cone over a suitably embedded
            $\P^3$, or $E$ is the Hirzebruch surface
            $F_1:=\P(\mathcal{O}_{\P^1}\oplus\mathcal{O}_{\P^1}(-1))$
            and $\overline{Y}$ \emph{is rigid};
            \item[$k=9$] then $E\cong \P^2$ and $\overline{Y}$ \emph{is
            rigid} (this case follows also by \cite{Schlessinger71}).
        \end{description}
    \end{description}
    \item On the other hand $E$ is a surface embedded in the smooth
    3--fold $Y$, which means that $E$ cannot admit
    non-hypersurface singularities. This fact gives significative
    constraints on the non--normal case implying that:
    \begin{itemize}
        \item \emph{if $E$ is non--normal then it is a suitable projection of a Hirzebruch surface
        $F_a := \P(\mathcal{O}_{\P^1}\oplus\mathcal{O}_{\P^1}(-a))$ having $E^3=\deg
        E=7$} (\cite{Gross97a}, Theorem 5.2). \emph{In this particular case there exists a
        smoothing $\widetilde{Y}$ of $\overline{Y}$}
        (\cite{Gross97a}, Lemma 5.6).
    \end{itemize}
\end{enumerate}

\begin{conclusion}
If \  $T(Y,\overline{Y},\widetilde{Y})$ is a type II transition
then $Y$ is the blow up of $\overline{Y}$ at the singular point
and the exceptional divisor $E$ is either a rational, normal, del
Pezzo surface of degree $k\leq 8$ or a non--normal del Pezzo
surface of degree 7. In the first case if
\begin{description}
    \item[$k\leq 3$] then $\overline{Y}$ has a hypersurface
    singularity,
    \item[$k=4$] then $\overline{Y}$ has a complete intersection
    singularity,
    \item[$k=5$] then $\overline{Y}$ is locally a Pfaffian
    subscheme,
    \item[$6\leq k\leq 8$] then $\overline{Y}$ is locally a cone
    over $E$ admitting a toric structure.
\end{description}
In particular $E$ can never be isomorphic to either $\P^2$ or
$F_1$ \emph{(\cite{Gross97a} Theorem 5.8)}.
\end{conclusion}

\subsubsection{Transitions of type III}

$\phi$ is the contraction of an irreducible divisor $E$ down to a
smooth curve $C\subset\overline{Y}$. Then:
\begin{enumerate}
    \item $C=\Sing (\overline{Y})$ is a smooth curve of \emph{canonical
    singularities} of $\overline{Y}$; apply Theorem 2.2 of
    \cite{Reid80} to conclude that $C$ is entirely composed of cDV
    singular points since $E$ is essentially the only possible
    exceptional divisor of a crepant resolution of $\overline{Y}$ and
    it gives a 1--dimensional fibration over $C$;
    \item the restriction $\phi|_E : E\rightarrow C$ exhibit $E$
    like a conic bundle over $C$, whose fibre is either a smooth
    conic, a union of two lines meeting at a point, or a double
    line; in particular if the general fibre is smooth then $E$ is
    normal (\cite{Wilson92} Theorem 2.2, \cite{Wilson93});
    \item let $\widehat{E}$ be the normalization of $E$ and
    $f:\widehat{E}\rightarrow Y$ the induced map; saying $Def(f)$
    the deformations space of $f$ like in \cite{Ran89} and
    $Def(Y)$ the Kuranishi space of $Y$, there is a natural map
    \begin{equation*}
    Def(f)\longrightarrow Def(Y)\ ;
    \end{equation*}
    then: \emph{the genus of $C$ is less or equal to the codimension of $\im (Def(f)\rightarrow
    Def(Y))$} (\cite{Gross97b}, Proposition 1.2);
    \item by the previous step: \emph{if $g(C)\geq 1$ then there exists a smoothing $\widetilde{Y}$ of
    $\overline{Y}$} (\cite{Gross97b}, Theorem 1.3); in fact there exists a deformation
    $\mathcal{Y}\rightarrow\Delta$
    of $Y$ such that the exceptional divisor $E$ do not deform to
    general $\mathcal{Y}_t, t\in \Delta$ since
    \begin{equation*}
    \codim \left(\im (Def(f)\rightarrow
    Def(Y))\right)\geq 1 \ ;
    \end{equation*}
    the contraction $\phi$ yields a contraction
    $\mathcal{Y}\rightarrow\mathcal{\overline{Y}}$ where
    $\mathcal{\overline{Y}}\rightarrow\Delta$ is the deformation
    induced by $\mathcal{Y}$ via the natural map $Def
    (Y)\rightarrow Def(\overline{Y})$, which exists by
    \cite{Kollar-Mori92}, Proposition 11.4; for general $t\in \Delta$ the contraction
    $\mathcal{Y}_t\rightarrow\mathcal{\overline{Y}}_t$ is then of
    type I; by \ref{tipo I} there is a smoothing
    $\mathcal{\widetilde{Y}}_t$ of $\mathcal{\overline{Y}}_t$
    except when $\Sing (\mathcal{\overline{Y}}_t)$ is composed by a unique
    ordinary double point; some more technical consideration shows
    that the latter does not occur for general $t$;
    \item it remains to understand what happens when $g(C)=0$
    i.e. $C\cong \P^1$; the goal is \emph{to construct a deformation
    $\mathcal{\overline{Y}}\rightarrow\Delta$ of
    $\overline{Y}$ such that the image of the induced map $\Delta\rightarrow
    Def(\overline{Y})$ is not contained in} $\im \left( Def(Y)\rightarrow
    Def(\overline{Y})\right)$; if such a deformation exists then
    $\mathcal{\overline{Y}}_t$ has $\Q$--factorial terminal
    singularities for general $t \in \Delta$ (\cite{Gross97b}, Lemma
    1.6) and by results of Y.~Namikawa and J.~Steenbrink \cite{Namikawa-Steenbrink95} it suffices
    to guarantee the existence of a smoothing
    $\mathcal{\widetilde{Y}}_t$ of $\mathcal{\overline{Y}}_t$; to
    show the existence of the deformation $\mathcal{\overline{Y}}$ needs a
    careful analysis of the structure of $Def(\overline{Y})$ and
    of the differential of the map $Def(Y)\rightarrow
    Def(\overline{Y})$:
    \begin{itemize}
        \item \emph{if $E^3\leq 6$ the cokernel of the above differential
        has dimension $\geq 2$; $Def(\overline{Y})$ is smooth when $E^3 \leq 5$;
        if $E^3=6$ then $Def(\overline{Y})$ may not be smooth but it is set--theoretically
        defined by at most one equation in a neighborhood of the origin of its tangent space;
        then the desired deformation $\mathcal{\overline{Y}}\rightarrow\Delta$ exists for $E^3\leq
        6$} (\cite{Gross97b}, Theorem 1.7).
    \end{itemize}
\end{enumerate}

\begin{conclusion}
If\quad  $T(Y,\overline{Y},\widetilde{Y})$ is a transition of type
III then the associated contraction $\phi$ fibers its exceptional
divisor $E$ as a conic bundle over the smooth curve $C=\Sing
(\overline{Y})$. Moreover $C$ is a locus of cDV singularities of
$\overline{Y}$ and either $g(C)\geq 1$ or $g(C)=0$ and $\deg E
\leq 6$.

\noindent In particular $\phi$ cannot fibre $E$ as a conic bundle
of degree 7 or 8 over $\P^1$ \emph{(\cite{Gross97b}, Theorem
0.4)}.
\end{conclusion}

\subsubsection{What about a general transition?}

The case of a general geometric transition
$T(Y,\overline{Y},\widetilde{Y})$ is much more complicated than
the case of a primitive one, essentially for two reasons:
\begin{itemize}
    \item the geometry of the exceptional locus $E$ can be very
    intricate,
    \item $\overline{Y}$ can then assume very general canonical
    singularities so that $Def(\overline{Y})$ can be very singular
    and the deformation theory of $\overline{Y}$ very
    complicated.
\end{itemize}
Some partial result can be obtained from Theorem \ref{ng-thm} or a
generalization of it in the case of non--complete intersection
singularities (see \cite{Gross97a}, definition 4.2 and Theorem
4.3): anyway $\overline{Y}$ is assumed to be $\Q$--factorial and
only admitting (a particular kind) of isolated singularities.

Moreover let us conclude by observing that, given a geometric
transition
\begin{equation*}
    \xymatrix@1{Y\ar@/_1pc/ @{.>}[rr]_T\ar[r]^{\phi}&
                \overline{Y}\ar@{<~>}[r]&\widetilde{Y}}
\end{equation*}
even the decomposition of $\phi$ in primitive contractions can be
non--invariant with respect to deformations of $Y$. In fact
\emph{if $\phi$ factors through a transition of type III then the
\ka cone may jump under deformation} (\cite{Wilson92},
\cite{Wilson93} main theorem, \cite{Namikawa94} Theorem C).

\section{The \cy web}

\subsection{Reid's fantasy}

An immediate consequence of Theorem \ref{cambio omologico} is
that, starting from a given \cy 3--fold $Y$, a conifold transition
produce a \emph{topologically distinct} \cy 3--fold
$\widetilde{Y}$. Actually there are plenty of topologically
distinct well known examples of \cy 3--folds and this fact seems
to definitely exclude the possibility of any kind of
``irreducibility" for any more or less defined concept of moduli
space of \cy 3--folds.

\noindent This is something new with respect to what happens in
the lower dimensional cases of elliptic curves and K3 surfaces.

\begin{description}
    \item[Elliptic curves] Any 1--dimensional compact complex
    manifold with $K_C\equiv 0$ is biholomorphic to an algebraic
    smooth plane cubic curve, i.e. to a complex torus, and viceversa.
    In particular their complex moduli space is the moduli space
    of complex structures over the topological torus $S^1 \times
    S^1$. Such a moduli space is algebraic,
    smooth and irreducible (the well known \emph{modular curve}).
    \item[K3 Surfaces](See \cite{Beauville78} and \cite{BPvdV84})
    The following facts were known to F.~Enriques \cite{Enriques46}:
    \begin{itemize}
        \item $\forall g \geq 3$ there exists a K3 surface of
        degree $2g-2$ in $\P^g$; hence its sectional genus is $g$;
        \item $\forall g \geq 3$ we can obtain a space $\mathcal{M}_g$
        of complex projective moduli of such
        surfaces, by imposing a polarization: $\mathcal{M}_{g}$ is an irreducible, analytic
        variety with $\dim_{\C}\mathcal{M}_{g}=19$;
        \item then the complex moduli space $\mathcal{M}^{alg}$ of algebraic K3
        surfaces is a \emph{reducible} analytic variety and it
        admits a countable number of irreducible components;
        \item there exist K3 surfaces belonging to more than one
        irreducible component of $\mathcal{M}^{alg}$; anyway if we
        restrict to K3's admitting $\Pic \cong \Z$ (they give the
        general element of any irreducible component) then they
        belong to only one irreducible component.
    \end{itemize}
    What could appear to F.~Enriques as a wildly reducible moduli space
    was explained by K.~Kodaira \cite{Kodaira64} as an analytic
    codimension 1 subvariety of a smooth, irreducible, analytic  variety
    $\mathcal{M}$. More precisely:
    \begin{itemize}
        \item there exist analytic non-algebraic K3 surfaces,
        \item the Kuranishi space of any analytic K3 surface is
        smooth and of dimension 20.
    \end{itemize}
    The latter suffices to construct a smooth, irreducible, analytic
    universal family of K3 surfaces: its base $\mathcal{M}$ is the complex
    analytic moduli space of K3 surfaces and $\dim_{\C}\mathcal{M}=20$.
    Moreover $\mathcal{M}^{alg}$ turns out to be a dense subset
    of $\mathcal{M}$.

    \noindent In other words \emph{the irreducibility of the moduli space of
    K3's is obtained by leaving the algebraic geometric category to work in
    the larger category of compact, \ka, analytic manifolds.} In fact \emph{any K3 surface is
    \ka} since all of them admit a canonical Ricci flat
    \ka--Einstein metric.
\end{description}

\noindent In \cite{Reid87} M.~Reid suggested that the right
approach to perceive some kind of irreducibility of a suitable
moduli space of \cy 3--folds could be similar to the case of K3
surfaces: one has to work in the right category. The key idea is
given by the following result of R.Friedman:

\begin{theorem}[\cite{Friedman86}, Corollary 4.7]\label{conifold vs. non-ka}
Let $\phi:Y\rightarrow\overline{Y}$
be a small contraction of a \cy 3--fold $Y$ to a normal 3--fold
$\overline{Y}$ such that $H^2(Y)$ is generated by the exceptional
locus $E$ of $\phi$ and $\Sing (\overline{Y})$ is composed by
$N\geq 2$ nodes. Then $\overline{Y}$ is smoothable and every
smoothing $\widetilde{Y}$ has $b_2(\widetilde{Y})=0$. Hence
$\overline{Y}$ can be smoothed only to \emph{non--\ka} compact
complex 3--folds.
\end{theorem}

\begin{corollary}There exist ``non--\ka \cy" 3--folds which can be
realized, by means of a conifold transition, starting from an
algebraic \cy 3--fold $Y$ as in Theorem \ref{conifold vs. non-ka}.
\end{corollary}

\begin{remark}
There is an evident contradiction in the words \emph{non--\ka \cy}
since in the definition \ref{cy-def} we assumed a projective
embedding for $Y$. Anyway their meaning should be evident as well
and probably the reader will forgive such an abuse of notation!
\end{remark}

A ``\cy" 3--fold with second Betti number equal to zero has
topological type completely determined by the third Betti number.
By results of C.~T.~C.~Wall \cite{Wall66} this suffices to guarantee
that it is diffeomorphic to a connected sum $\left(
S^3\times S^3\right) ^{\# r}$ of $r$ copies of the \emph{solid
hypertorus} $S^3\times S^3$. Introduce then the following:

\begin{assumptions}
\begin{enumerate}
    \item every projective \cy 3--fold $Y$ is birational to a \cy 3--fold
    $Y'$ such that $H^2(Y')$ is generated by rational
    curves; moreover if $\phi:Y'\rightarrow\overline{Y}$ is the
    morphism contacting all them, then $\overline{Y}$ is always smoothable;
    \item the moduli space $\mathcal{N}_r$ of complex structures
    on $\left( S^3\times S^3\right) ^{\# r}$ is
    \emph{irreducible}.
\end{enumerate}
\end{assumptions}

Then we get the famous:

\begin{conjecture}[the Reid's fantasy]\label{fantasia di Reid}
Up to some kind of inductive limit over $r$, the birational
classes of projective \cy 3--folds can be fitted together, by
means of geometric transitions, into one irreducible family
parameterized by the moduli space $\mathcal{N}$ of complex
structures over suitable connected sum of copies of solid
hypertori.
\end{conjecture}

In fact if $Y$ is a \cy 3--fold, by assumption (1) we can recover
a birational \cy 3--fold $Y'$ admitting a small contraction
morphism $\phi:Y'\rightarrow\overline{Y}$. Since $\phi$ is a
composition of a finite number of type I contractions, \ref{tipo
I} guarantees that $\overline{Y}$ admits at most a finite number
of isolated cDV singular points.

\noindent Then by Theorem \ref{ng-thm}, $\overline{Y}$ can be
deformed to a variety $\overline{Y}'$ admitting at worst nodes as
singularities. Recalling Theorem \ref{friedman-thm}, the second
part of assumption (1) implies that either $|\Sing \left(
\overline{Y}'\right) |\geq 2$ or $\overline{Y}'$ is smooth. In the
first case Theorem \ref{friedman-thm}, or equivalently Theorem
\ref{conifold vs. non-ka}, gives a smoothing $\widetilde{Y}$ of
$\overline{Y}'$. In the second case rename $\overline{Y}'$ as
$\widetilde{Y}$. In both cases $\widetilde{Y}$ is a non--\ka \cy
3--fold since $H^2(\widetilde{Y})=0$. Then it is diffeomorphic to
a connected sum of $r$ copies of solid hypertori, where $r$
\emph{depends on the topology of }$Y$. In fact, if in particular
we make the further assumption that $\Sing(\overline{Y})$ \emph{is
composed only by nodes} then the transition
$T(Y',\overline{Y},\widetilde{Y})$ is a conifold one and Theorem
\ref{cambio omologico} gives
\begin{equation*}
    r=b_3(\widetilde{Y})/2 = b_3(Y')/2 + c = b_3(Y')/2+N-k
\end{equation*}
Assumption (1) implies that $k=b_2(Y')$ and that $N\geq 2$. The
previous relation can be then rewritten as follows:
\begin{equation}\label{vincolo top}
    b_3(Y')-2b_2(Y')=2r- N \leq 2r - 2
\end{equation}
Since $Y$ and $Y'$ are birational, their Betti numbers
coincides\footnote{This is a famous result of V.~Batyrev
\cite{Batyrev99}, obtained by employing $p$--adic integration and
Weil conjectures. It seems that this result motivated M.~Kontsevich
to introduce the theory of \emph{motivic integration} in a
memorable lecture at Orsay \cite{Kontsevich95}, in which he proved
that two birational  \cy varieties even have isomorphic Hodge
structures. Actually, as explained by Batyrev in the introduction
of \cite{Batyrev99}, the 3--dimensional case, to which we are
interested here, can be deduced by an older result of Y.~Kawamata
\cite{Kawamata88}, since two birational minimal models of 3--folds
can be connected by a sequence of flops.}. Then (\ref{vincolo
top}) can be rewritten, in terms of the Euler--Poincar\'e
characteristic of $Y$, as follows
\begin{equation}\label{vincolo top-bis}
    r = 1 + \frac{N - \chi (Y)}{2}\geq 2 - \chi (Y)/2
\end{equation}
In conclusion, by means of a geometric transition, the birational
equivalence class of the \cy 3--fold $Y$ determines a complex
structure over $\left( S^3\times S^3\right) ^{\# r}$, given by
$\widetilde{Y}$ and represented by a point of $\mathcal{N}_r$, for
$r\gg 0$ according with (\ref{vincolo top-bis}). On the other
hand, results stated in \ref{cambio omologico generalizzato}
ensure that the previous argument applies, with slight
modifications, to any \cy 3--fold, without the assumption that $T$
is conifold.

\noindent The last step should be a sort of \emph{gluing} of all
the $\mathcal{N}_r$'s preserving irreducibility postulated by
assumption (2) (to use M.~Reid's words: ``let's ignore this as a
minor technical problem").

\begin{remark}\label{non-ka}
The key point of the Reid's fantasy is clearly the assumption (2):
very little is known about complex structures over solid hypertori
and very few techniques are available in dealing with compact
complex non--\ka manifolds!
\end{remark}

\begin{remark}
The geometric beauty of the Reid's fantasy \ref{fantasia di Reid}
is given also by the evident analogies with both the lower
dimensional cases of elliptic curves and of K3 surfaces. In fact
as in the last case, the irreducibility of the moduli space is
recovered by means of particular geometric transitions which
actually are the right tools \emph{to leave the compact, \ka
category to work into the larger category of compact, complex,
analytic manifolds}. On the other hand, as in the case of elliptic
curves, the moduli problem is reduced to parameterize complex
structures over a sort of ``generalized tori".
\end{remark}

\subsection{The ``vacuum degeneracy problem" in string theory}

The geometric transition's property of connecting topologically
distinct \cy 3--folds and in particular the restored concept of
 a possible irreducible moduli space due to the Reid's Conjecture
 \ref{fantasia di Reid} suggested most interesting applications in
 string theory.

 \noindent In fact \cy 3--folds play a fundamental role in
 10--dimensional string theories: locally 4 dimensions give rise
 to the usual Minkovsky space--time while the remaining 6 dimensions (the so called
 \emph{hidden dimensions} for their microscopic extension, of the same order as the
 Plank constant) are compactified to a geometric model
 which, essentially to preserve the required supersymmetry, turns out to
 be a \cy 3--fold.

 \noindent In spite of the fact that there are very few
 consistent 10--dimensional super--string theories, actually
 near--unique via dualities, the compactification process give
 rise to the problem of choosing the appropriate \cy model: on one hand there
 is not any prescription for making a precise choice and on the
 other hand there is a huge multitude of topologically distinct
 \cy 3--folds. Moreover the choice of two distinct \cy models is
 not ``a priori" equivalent from the physical point of view, since
 the second and the third Betti numbers (or better the Hodge numbers
 $h^{1,1}$ and $h^{2,1}$) of the \cy model are strictly related with the number of
 hypermultiplets and the number of vector multiplets,
 respectively, of the compactified physical theory.

 \noindent This is the so called \emph{vacuum degeneracy problem}
 in string theory.

 \noindent The ideas of Clemens and then of Friedman and Reid,
 leading to the formulation of the Reid's fantasy in 1987
 suggested to physicists like P.~Candelas, P.~S.~Green, T.~H\"{u}bsch
 and others that:
 \begin{itemize}
    \item \emph{\cy 3--folds could be, at least mathematically,
    connected each other by means of geometric (conifold) transitions}.
 \end{itemize}
This is the so called \emph{\cy web conjecture} described in many
insightful papers starting from 1988 (see \cite{CDLS88},
\cite{Green-Hubsch88let}, \cite{Green-Hubsch88}, \cite{CGH89},
\cite{CGH90}). A more precise version of this conjecture will be
given later following M.~Gross (see \ref{cong. di connessione}).

\noindent In the previous statement \emph{mathematically} means
that the geometric (or eventually the conifold) transition
connecting each other two \cy 3--fold is merely a geometrical
process: what about the \emph{physical transition} between the
physical theories involved?

\noindent A first answer was given, for what concerning a conifold
transition, in 1995 by A.~Strominger (see \cite{Strominger95} and
\cite{GMS95}). His explanation of how physical theories can pass
\emph{smoothly} through the conifold singularities of the moduli
space of \cy string vacua was inspired by techniques of N.~Seiberg
and E.~Witten \cite{Sieberg-Witten94}: the idea is that the
topological change is given by the \emph{condensation} of massive
black holes to massless ones.

\noindent In the following years some other geometric transition,
more general than the conifold one, have been physically
understood: see for example \cite{BKK95}, \cite{KMP96},
\cite{BKKM97}.

\subsection{The connectedness conjecture}\label{cong. di connessione}

A mathematically refined version of the \cy web conjecture was
presented by M.~Gross in \cite{Gross97b}.

\noindent On the contrary of the K3 case for which an algebraic K3
surface can be smoothly deformed to a non--algebraic one, the
deformation of a projective \cy 3--fold, even singular, is still
projective. Since the hardest part of the Conjecture \ref{fantasia
di Reid} seems to be in dealing with non--\ka \cy 3--folds and in
finding non--algebraic contractions, as observed in Remark
\ref{non-ka}, one could skip this part by insisting on staying
within the projective category as follows.

\noindent One can think the nodes of the giant web predicted by
the web conjecture as consisting in \emph{deformation classes} of
\cy 3--folds. Two of such nodes, say $\mathcal{M}_1$ and
$\mathcal{M}_2$, are \emph{connected} by an arrow
$\mathcal{M}_1\rightarrow\mathcal{M}_2$ if \emph{the general
element of $\mathcal{M}_1$ is connected with a smooth element of
$\mathcal{M}_2$ by means of a geometric transition}, which means:
for the general element $Y$ of $\mathcal{M}_1$ there exists a
birational contraction to a normal 3--fold
$\phi:Y\rightarrow\overline{Y}$ and a flat local family
$\overline{\mathcal{Y}}\rightarrow\Delta$ whose central fibre is
$\overline{\mathcal{Y}}_0\cong \overline{Y}$ and such that
$\overline{\mathcal{Y}}_t$ is a smooth element of $\mathcal{M}_2$
for general $t\in\Delta$.

\begin{example}[See also \cite{Gross97b}]\label{esempio-web} Let
\begin{itemize}
    \item $\mathcal{M}_Q$ be the moduli space of smooth quintic 3--folds in $\P
    ^4$,
    \item $\mathcal{M}_D$ be the moduli space of double solids (i.e. double
    covers of $\P ^3$) branching along a smooth octic surface of $\P ^3$,
    \item $\mathcal{M}_T$ be the moduli space of smooth blow--up's of
    quintic 3--folds having a triple point.
\end{itemize}
Let $Z$ be a general element in $\mathcal{M}_T$ and
$\phi:Z\rightarrow\overline{Y}$ be the contraction of the
exceptional divisor of $Z$. Then $\overline{Y}$ is a quintic
3--fold in $\P ^4$ with a triple point. Since $\overline{Y}$ can
be smoothed to a quintic 3--fold we have
\begin{equation}\label{M_T>M_Q}
    \mathcal{M}_T\longrightarrow\mathcal{M}_Q
\end{equation}
by means of a primitive transition of type II.

\noindent On the other hand if we project $\overline{Y}$ from the
triple point $p_o$ we get a rational morphism
\begin{equation*}
    \psi: \overline{Y}\dashrightarrow\P ^3
\end{equation*}
\begin{proposition}\label{8ds}
The previous rational morphism $\psi$ can be lifted to the blow up
$Z$ giving rise to a generically finite morphism $\widehat{\psi}:
Z\rightarrow \P ^3$. More precisely $\widehat{\psi}$ is 2:1 except
over 60 points $\{ p_i\}$ for which $\widehat{\psi}^{-1}(p_i)\cong
\P ^1$. Consider the Stein factorization
$\widehat{\psi}=f\circ\varphi$. Then we get the following
commutative diagram
\begin{equation}\label{diagrammaStein}
    \xymatrix{Z\ar[r]^{\varphi}\ar[dr]^{\widehat{\psi}}\ar[d]^{\phi}& \overline{X}\ar[d]^{f}\\
              \overline{Y}\ar@{-->}[r]^{\psi}& \P ^3}
\end{equation}
where $\varphi$ is the birational contraction of all of the 60 $\P
^1$'s and $f$ gives to $\overline{X}$ the structure of a double
solid branched along a singular octic surface $S\subset\P ^3$.
\end{proposition}
\noindent Since $\overline{X}$ can immediately be smoothed by
smoothing the branching locus $S\subset\P ^3$ it is possible to
write
\begin{equation}\label{M_T>M_D}
    \mathcal{M}_T\longrightarrow\mathcal{M}_D
\end{equation}

\begin{itemize}
    \item \emph{Therefore the deformation families
$\mathcal{M}_Q,\mathcal{M}_T,\mathcal{M}_D$ are nodes of the
following \emph{connected graph} obtained by composing
(\ref{M_T>M_Q}) and (\ref{M_T>M_D})}:
\begin{equation}\label{3web}
    \xymatrix{&\mathcal{M}_T\ar[dr]\ar[dl]&\\
              \mathcal{M}_Q&&\mathcal{M}_D}
\end{equation}
\end{itemize}

\begin{proof}[Proof of Proposition \ref{8ds}]
The rational morphism $\psi$ is defined as follows
\begin{equation*}
    \forall p \in \overline{Y}\setminus\{ p_o\}\quad \psi (p):=
    l(p_o,p)\in \mathcal{L}_{p_o}:=\left\{\text{lines $l$ of $\P ^4$ through $p_o$}\right\}\cong \P
    ^3
\end{equation*}
where $l(p_o,p)$ is the line connecting $p$ and $p_o$. Since the
domain $\overline{Y}\setminus\{ p_o\}$ of $\psi$ coincides with
the locus of smooth points of $\overline{Y}$, $\psi$ can be
naturally lifted to a well defined morphism
$\widehat{\psi}:Z\rightarrow \P ^3$ by setting
\begin{eqnarray*}
  \forall q\in Z\setminus E\quad \widehat{\psi}(q) &=& (\psi\circ\phi) (q) \\
  \forall q\in E\quad \widehat{\psi}(q) &=& l_q
\end{eqnarray*}
where $E$ is the exceptional locus of the blow up $\phi$ and $l_q$
is the tangent line to $\overline{Y}$ in $p_o$ determined by the
tangent direction represented by $q\in E$. The morphism
$\widehat{\psi}$ is clearly generically 2:1 and the image in
$\mathcal{L}_{p_o}\cong\P ^3$ of the branching locus is given by
\begin{equation*}
    S:=\left\{\text{lines $l(p_o,p)$ which are tangent to $\overline{Y}$ in
    $p$}\right\}\subset\P ^3
\end{equation*}
\emph{$S$ is a surface of degree 8}. In fact locally the triple
point $p_o$ can be assumed to be the origin of an affine subset
$\C ^4$ of $\P ^4$. The local equation of $\overline{Y}$ is then
given by $F_5+F_4+F_3=0$ where $F_d=F_d(x,y,z,w)$ is a generic
homogeneous polynomial of degree $d$. If $p=(x_p,y_p,z_p,w_p)$
then $l(p_o,p)$ is parameterized by
\begin{equation*}
    x=x_pt\ ,\ y=y_pt\ ,\ z=z_pt\ ,\ w=w_pt
\end{equation*}
Therefore $l(p_o,p)\in S$ if and only if
\begin{equation*}
   (F_5+F_4+F_3)|_{l(p_o,p)}=t^3(at^2+bt+c)
\end{equation*}
where $a,b,c$ are homogeneous polynomials in $x_p,y_p,z_p,w_p$ of
degree 5,4,3, respectively, satisfying the further tangency
condition
\begin{equation}\label{delta=0}
    b^2-4ac=0\ .
\end{equation}
The latter gives a degree 8 homogeneous equation in $\P
^3(x_p,y_p,z_p,w_p)\cong \mathcal{L}_{p_o}$.

\noindent Observe that the 60 points $\{ p_i\}$ described in
$\mathcal{L}_{p_o}$ by $a=b=c=0$ are the images via $\psi$ of the
lines contained in $\overline{Y}$. Hence
$\widehat{\psi}^{-1}(p_i)\cong \P ^1$ while $\widehat{\psi}$ is
2:1 over $\P ^3 \setminus \{ p_i\}$. The Stein factorization
$\widehat{\psi}=f\circ\varphi$ is then the composition of the
birational morphism $\varphi$ contracting all of those $\P ^1$'s
and of the 2:1 morphism $f$ onto $\P ^3$.

\noindent The situation is then described by the commutative
diagram (\ref{diagrammaStein}) where $\overline{X}$ is a double
covering of $\P ^3$ branched along the surface $S$. Since equation
(\ref{delta=0}) of $S$ gives $\Sing (S) = \{ a=b=c=0\}$,
$\overline{X}$ admits the 60 isolated singularities given by the
images by $\varphi$ of the contracted $\P ^1$'s. The smoothing of
$\overline{X}$ is then given by the double solid branched along
the generic surface of degree 8 in $\P ^3$.
\end{proof}
\end{example}

Let us come back to the connected graph (\ref{3web}). Then the
question is: can that graph be enlarged to a very bigger graph
connecting deformation classes of all simply connected \cy
3--folds?

\begin{conjecture}[of Connectedness]\label{connessione}
The graph of simply connected \cy 3--folds is connected.
\end{conjecture}

Evidences for such a conjecture were firstly given in
\cite{Green-Hubsch88}, where the moduli spaces of some \cy
3--folds, which are complete intersections in products of
projective spaces, were connected each other.

\noindent Most significative evidences are given in \cite{CGGK96}
where a general procedure for connecting up \cy 3--folds which are
complete intersections in some toric variety, is described. Such a
procedure was developed starting from an original idea of
D.~Morrison and works by \emph{intersecting} the combinatorial
toric data (i.e. \emph{reflexive polytopes}) of two given \cy
3--folds, to produce a further \cy 3--fold (\emph{if} the so
obtained toric data give rise to a reflexive polytope too!). The
latter \cy is then connected to the previous two, by means of
geometric transitions. By direct computer search, the authors
checked that the procedure described allows to settle all known
examples of \cy hypersurfaces in weighted $\P ^4$ (7555 \cy
3--folds) into a big connected graph. This result was actually
already known to P.~Candelas and collaborators, but the new fact is
that the third \cy 3--fold, obtained \emph{by intersecting} the
toric data of two given \cy weighted hypersurfaces, is not, in
general, a weighted hypersurface but rather a complete
intersection in a more general toric variety. Which is: \emph{the
graph connecting up all the 7555 \cy weighted hypersurfaces
extends to englobe many complete intersections in more general
toric varieties}.

\noindent Let us remark that, in general, the geometric
transitions involved in the procedure described above are not
conifold. Hence such a big graph produces a \emph{mathematical}
link between deformation classes of \cy 3--folds, leaving open the
problem of a satisfying physical understanding of the induced
connection between string vacua.

\section{Mirror symmetry and transitions: the reverse transition}

A natural question arises from the previous connectedness
Conjecture \ref{connessione}:
\begin{itemize}
    \item \emph{is such a conjecture consistent with already known
    ``connecting processes" between \cy string vacua suggested by
    physical dualities like e.g. \emph{mirror symmetry}?}
\end{itemize}
In a sense, a positive answer to this question represents a
further evidence supporting the stated conjecture.

\subsection{Mirror symmetry conjecture: some mathematical statements}

A description of physical origin and meaning of mirror symmetry
conjecture is outside the scope of this paper. In the following we
will simply state some (minimal) mathematical consequences useful
to understand the role of geometric transition in this context.
The reader interested in a deeper understanding of the topic
should consult the extensive monographs \cite{Voisin96},
\cite{Cox-Katz99} and the recent \cite{HKKPTVVZ03}.

\begin{conjecture}[Infinitesimal Mirror Symmetry]\label{ms-infinitesimale}
Let $Y$ be a \cy variety. Then there exists a \cy variety
$Y^{\circ}$ and isomorphisms of complex vector spaces
\begin{equation}\label{inf.mir.sym}
    \forall\ 0\leq p, q \leq \dim Y\quad
    \xymatrix@1{\mu_{p,q}:H^p \left(\Omega^q_{Y}\right)\ar[r]^-{\cong}&H^p
    \left(\Omega^{n-q}_{Y^{\circ}}\right)}
\end{equation}
inducing a mirror reversing identification on the Hodge diamonds
of $Y$ and $Y^{\circ}$.
\end{conjecture}

\begin{remark}
Since $\mathcal{K}_Y \cong \mathcal{O}_Y$ we get canonical
isomorphisms
\begin{equation*}
     H^p \left( \Omega^{n-q}_Y\right)\cong H^p \left( \bigwedge^q \mathcal{T}_Y\right)
\end{equation*}
and the same for $Y^{\circ}$. Their composition with isomorphisms
$\mu_{p,q}$ in (\ref{inf.mir.sym}) give rise to the following
isomorphisms
\begin{eqnarray*}
    \forall\  0\leq p,q \leq \dim Y\quad
  &&\xymatrix@1{\mu_{p,q}': H^p \left( \bigwedge^q \mathcal{T}_Y\right)\ar[r]^-{\cong}&
                                                H^p \left(\Omega^q_{Y^{\circ}}\right)}    \\
  &&\xymatrix@1{\mu_{p,q}'': H^p \left(\Omega^q_{Y}\right)\ar[r]^-{\cong}&
                                     H^p \left( \bigwedge^q \mathcal{T}_{Y^{\circ}}\right)}
\end{eqnarray*}
and commutative diagrams
\begin{equation}\label{ms-diagramma}
    \xymatrix{H^p \left( \bigwedge^{n-q} \mathcal{T}_Y\right)\ar[r]^-{\mu_{p,n-q}'}\ar[d]^{\cong}&
                            H^p \left(\Omega^{n-q}_{Y^{\circ}}\right)\ar[d]^{\cong}\\
              H^p \left(\Omega^q_{Y}\right)\ar[r]^-{\mu_{p,q}''}\ar[ru]^{\mu_{p,q}}&
                            H^p \left( \bigwedge^q
                            \mathcal{T}_{Y^{\circ}}\right)}\ .
\end{equation}
In particular if $p=1=q$ then
\begin{eqnarray}\label{mirror-diff.}
  &&\xymatrix@1{\mu':=\mu_{1,1}': H^1 \left( \mathcal{T}_Y\right)\ar[r]^-{\cong}&
                                                H^1 \left(\Omega_{Y^{\circ}}\right)}    \\
  &&\xymatrix@1{\mu'':=\mu_{1,1}'': H^1 \left(\Omega_{Y}\right)\ar[r]^-{\cong}&
                                     H^1 \left( \mathcal{T}_{Y^{\circ}}\right)}
\end{eqnarray}
\end{remark}

\subsubsection{The \cy moduli space}\label{cy-moduli}

To give a \cy variety $Y$ means in particular to fix a
\emph{triple} $(Y,J,h)$ of a \emph{compact manifold} $Y$, a
\emph{complex structure} $J$ on $Y$ and a \emph{hermitian metric}
$h$ on $Y$ whose real part gives a \emph{Ricci flat} riemannian
metric, and whose imaginary part gives a closed $(1,1)$--form
$\omega:=-1/2\im h$ (i.e. a \emph{\ka form}) which is
\emph{positive}.

\noindent Think the \emph{complex moduli space}
$\mathcal{M}_Y^{\C}$ of $(Y,J,h)$ as the space parameterizing all
the deformations of the complex structure $J$ over $Y$ up to
biholomorphisms. The Bogomolov--Tian--Todorov theorem asserts that
\emph{locally $\mathcal{M}_Y^{\C}$ is smooth} (see
\cite{Bogomolov78}, \cite{Tian87}, \cite{Todorov89} and also
\cite{Ran92} for a more recent and algebraic proof). Then:
\begin{itemize}
    \item $H^1 \left( \mathcal{T}_Y\right)$ can be canonically
identified with the \emph{tangent space to} $\mathcal{M}_Y^{\C}$
at the fixed complex structure $J$.
\end{itemize}

\noindent On the other hand the Yau theorem solving the Calabi
conjecture (see \cite{Calabi56} and \cite{Yau78}) ensures that,
\emph{for any positive \ka form $\omega$ such that $[\omega]\in
H^2(Y,\R)\cap H^1(\Omega_Y)$, there exists a unique Ricci flat
metric whose associated $(1,1)$--form is cohomologous to
$\omega$.} Then Definition \ref{cono di kahler} and Corollary
\ref{kahler-interpretazione} imply that all the possible
deformations of the Ricci flat, \ka metric $h$ on $Y$ are
parameterized by the \ka cone $\mathcal{K}(Y)$. For this reason
$H^1 \left(\Omega_{Y}\right)$ can be thought as a
\emph{complexification of the tangent space to the \ka moduli
space of $Y$}.

\noindent Moreover one can give a more natural meaning to $H^1
\left(\Omega_{Y}\right)$ by constructing a \emph{complexified \ka
moduli space} as follows.

\noindent First of all observe that the mathematical datum of a
given \cy 3--fold $Y$, which is actually a triple $(Y,J,h)$ with
$\dim Y=3$, do not completely characterize the physical string
theory compactified to $Y$. To do this an extra--datum, called
\emph{the $B$--field}, is needed. Physically it is a
characteristic parameter of the string action. Mathematically it
is represented by the choice of a lateral class $\beta$ in the
quotient $H^2(Y,\R)/H^2(Y,\Z)$. One can then look at the complex
class
\begin{equation*}
    \chi:= \beta + i \omega = \beta - i/2 \im h \in
    H^2(Y,\C)/H^2(Y,\Z)
\end{equation*}
where $H^2(Y,\Z)$ acts naturally by inclusion, hence it acts only
on the real part of a class in $H^2(Y,\C)$, as desired. The class
$\chi$ is called the \emph{complexified \ka class} of the \cy
3--fold $Y$. It can be actually thought as a \emph{polarization}
over the complex variety $(Y,J)$ whose possible deformations are
then parameterized by the \emph{complexified \ka space}
\begin{equation*}
    \mathcal{K}_{\C}(Y):=\left\{ \chi\in H^2(Y,\C)|\im \chi \in
    \mathcal{K}(Y)\right\}/H^2(Y,\Z)\ .
\end{equation*}
The \emph{complexified \ka moduli space} $\mathcal{M}_Y^K$ is then
given by $\mathcal{K}_{\C}(Y)$ up to the action of the
automorphisms group $\Aut (Y)$ and
\begin{itemize}
    \item $H^1\left( Y,\Omega_Y\right)$ is the \emph{tangent space
    to} $\mathcal{M}_Y^K$ at the fixed complexified \ka class
    $\chi$,
    \item the \emph{\cy moduli space} of a \cy variety $Y$ is then
    the total space of a fibration
    \begin{equation}\label{fibrazione-moduli}
    \mathcal{M}_Y\longrightarrow\mathcal{M}_Y^{\C}
    \end{equation}
    whose fibre over the isomorphism class in $\mathcal{M}_Y^{\C}$
    represented by $(Y,J)$ is given by the complexified \ka
    moduli space $\mathcal{M}_Y^K$.
\end{itemize}
In particular, if $\dim Y = 3$, P.~M.~H.~Wilson proved that, outside
of a countable union of closed subsets of $\mathcal{M}_Y^{\C}$,
the \ka cone do not varies with the complex structure $J$ (see
\cite{Wilson92}, \cite{Wilson93}). That's enough to conclude that:
\begin{itemize}
    \item if $\dim Y = 3$ \emph{the fibration (\ref{fibrazione-moduli})
    is generically locally trivial}, which means that if $J$ is the
    complex structure of a sufficiently general \cy 3--fold $Y$
    then there exists a Zariski open subset
    $U\subset\mathcal{M}_Y^{\C}$ containing the class
    represented by $(Y,J)$ and such that $\mathcal{M}_Y|_U\cong
    U\times\mathcal{M}_Y^K$.
\end{itemize}
The Conjecture \ref{ms-infinitesimale} can then be understood as
the \emph{differential version} of the following one.

\begin{conjecture}[Local Mirror Symmetry for \cy 3--fold]\label{ms-locale}
Let $(Y,\chi)$ be the \emph{polarized couple} given by a general
\cy 3--fold $Y = (Y,J,h)$ and a complexified \ka class
$\chi\in\mathcal{K}_{\C}(Y)$ such that $\im \chi = -1/2\im h$.
Then there exist:
\begin{enumerate}
    \item a \emph{mirror polarized couple} $(Y^{\circ},\chi^{\circ})$,
    where $Y^{\circ}=(Y^{\circ},J^{\circ},h^{\circ})$ is a sufficiently general \cy 3--fold and
    $\chi^{\circ}\in\mathcal{K}_{\C}(Y^{\circ})$ is such that $\im \chi^{\circ} = -1/2\im h^{\circ}$,
    \item two open subsets $\mathcal{U}\subset\mathcal{M}_Y\ ,\
    \mathcal{U}^{\circ}\subset\mathcal{M}_{Y^{\circ}}$ containing
    the isomorphisms classes represented by $(Y,\chi)$ and
    $(Y^{\circ},\chi^{\circ})$, respectively; notice that they
    inherits the local product structure of $\mathcal{M}_Y$ and
    $\mathcal{M}_{Y^{\circ}}$ i.e.
    \begin{equation*}
    \mathcal{U}\cong U_{\C}\times U_K \quad,\quad
    \mathcal{U}^{\circ}\cong U^{\circ}_{\C}\times U_K^{\circ}
    \end{equation*}
    \item a biholomorphism
    $m:\mathcal{U}\rightarrow\mathcal{U}^{\circ}$, called \emph{local mirror
    map}, reversing the product structures, which is
    \begin{equation*}
    m\left( U_{\C}\right)=U_K^{\circ}\quad,\quad m\left( U_{K}\right)=U_{\C}^{\circ}
    \end{equation*}
    whose differential gives maps $\mu '$ and $\mu ''$ in
    (\ref{mirror-diff.}), i.e.
    \begin{equation*}
    d_{(J,\chi)}(m) = \mu '\times \mu ''
    \end{equation*}
\end{enumerate}
\end{conjecture}

\begin{remark}[Mirror partners of rigid \cy varieties]\label{rigid-cy}
Let $Y$ be a \emph{rigid} \cy variety i.e. $Y$ do not admits
complex deformations and
\begin{equation}\label{rigid-numbers}
    h^1(\mathcal{T}_Y)=h^{2,1}(Y)=0
\end{equation}
Assume $Y^{\circ}$ to be a mirror partner of $Y$. Then Conjecture
\ref{ms-infinitesimale} gives
\begin{equation}\label{mirror-rigid}
    h^{1,1}(Y^{\circ})=h^{2,1}(Y)=0
\end{equation}
which implies that $Y^{\circ}$ \emph{cannot be a \ka variety: in
particular $Y^{\circ}$ is not a \cy variety}.

\noindent Since rigid \cy 3--folds exist (the first examples were
constructed in 1986 by C.~Schoen in \cite{Schoen86}) this fact
introduces a counterexample to both the stated mirror symmetry
conjectures \ref{ms-infinitesimale} and \ref{ms-locale}.

\noindent From the mathematical point of view, such a
contradiction could be resolved by assuming mirror symmetry to
involve some \emph{non--\ka \cy} variety too (recall Remark
\ref{non-ka}): but which of them?

\noindent Anyway, from the physical point of view, it is
completely unclear which kind of string theory can be compactified
to a non--\ka \cy 3--fold: so what is the mirror dual of a string
theory compactified to a rigid \cy vacuum?
\end{remark}

\subsection{The reverse transition}\label{transizione reverse}

Consider a transition $T(Y,\overline{Y},\widetilde{Y})$ and let
$Y^{\circ}$ and $\widetilde{Y}^{\circ}$ be mirror partners of $Y$
and $\widetilde{Y}$, respectively:
\begin{equation}\label{MS+transition}
    \xymatrix{Y\ar@/^1pc/ @{.>}[rr]^T\ar[r]_{\phi}\ar@{<.>}[d]^{M.S.}
              &\overline{Y}\ar@{<~>}[r]&\widetilde{Y}\ar@{<.>}[d]^{M.S.}\\
               Y^{\circ}&&\widetilde{Y}^{\circ}}
\end{equation}
Recall that mirror symmetry exchange complex moduli with \ka
moduli. On the other hand, if $T$ is a conifold transition, point
(3) of Theorem \ref{cambio omologico} and Remark \ref{cpx<->ka}
allow to conclude that the topologies of $Y^{\circ}$ and
$\widetilde{Y}^{\circ}$ are compatible with a (reverse) conifold
transition
$T^{\circ}(\widetilde{Y}^{\circ},\overline{Y}^{\circ},Y^{\circ})$
which would complete diagram (\ref{MS+transition}) as follows
\begin{equation}\label{reverse}
    \xymatrix{Y\ar@/^1pc/ @{.>}[rr]^T\ar[r]_{\phi}\ar@{<.>}[d]^{M.S.}
              &\overline{Y}\ar@{<~>}[r]&\widetilde{Y}\ar@{<.>}[d]^{M.S.}\\
               Y^{\circ}\ar@{<~>}[r]&\overline{Y}^{\circ}&\widetilde{Y}^{\circ}
               \ar@/^1pc/@{.>}[ll]^{T^{\circ}}\ar[l]_-{\phi ^{\circ}}}
\end{equation}
Notice that the \emph{reverse conifold transition} $T^{\circ}$
would have the same parameters $N,k,c$ as $T$ whose role is now
reversed. Precisely
\begin{itemize}
    \item $\Sing \left(\overline{Y} ^{\circ}\right)$ would be composed by $N$ ordinary double
    points, just like $\Sing (\overline{Y})$,
    \item the exceptional locus of the birational contraction $\phi
    ^{\circ}$ would be composed by $N$ rational curves whose homology
    classes span a $c$--dimensional subspace of
    $H_2(\widetilde{Y}^{\circ})$,
    \item the vanishing locus of the smoothing $Y^{\circ}$ would
    be given by $N$ 3--spheres whose homology classes span a
    $k$--dimensional subspace of $H_3(Y^{\circ})$.
\end{itemize}
A similar picture naturally suggested that a diagram like
(\ref{reverse}) could be established for every geometric
transition $T$, leading to the following conjecture, probably due
to D.~Morrison.

\begin{conjecture}[of Reverse Transition, see \cite{Morrison99}, \cite{GMS95},
\cite{CGGK96} and \cite{Lynker-Schimmrigk}]\label{cong.reverse}
Let $T(Y,\overline{Y},\widetilde{Y})$ be a geometric transition
and let $Y^{\circ}$ and $\widetilde{Y}^{\circ}$ be mirror partners
of $Y$ and $\widetilde{Y}$, respectively. Then mirror partners are
linked by a \emph{reverse geometric transition}
$T^{\circ}(\widetilde{Y}^{\circ},\overline{Y}^{\circ},Y^{\circ})$
like in diagram (\ref{reverse}).
\end{conjecture}

In \cite{Morrison99} D.~Morrison supported such a conjecture with
an example employing the Greene--Plesser construction
\cite{Greene-Plesser90} to produce mirror partners of the
geometric transition linking a desingularization of an octic
weighted hypersurface of $\P (1,1,2,2,2)$ with the generic
complete intersection of bi--degree $(2,4)$ in $\P ^5$.

\noindent Further evidences were given in \cite{BC-FKvS98} where
the reverse transition of  a conifold transition, linking a
complete intersection in a Grassmannian with a complete
intersection in a Fano toric variety, is produced: in particular
the reverse transition is still conifold. This fact suggests to
specialize Conjecture \ref{cong.reverse} as follows.

\begin{conjecture}\label{cong.reverse-conifold} Let
$T(Y,\overline{Y},\widetilde{Y})$ be a conifold transition. Then
there exist mirror part\-ners $Y^{\circ}$ and
$\widetilde{Y}^{\circ}$ of $Y$ and $\widetilde{Y}$ and a re\-verse
tran\-si\-tion
$T^{\circ}(\widetilde{Y}^{\circ},\overline{Y}^{\circ},Y^{\circ})$
which is still conifold.
\end{conjecture}

Such a conjecture seems to be natural when we look at the role
played by parameters $N,k,c$. Anyway in \cite{KMP96} examples of
geometric \emph{non--conifold} transitions
$T(Y,\overline{Y},\widetilde{Y})$, which can be deformed to
conifold transitions, are produced. More precisely the birational
contraction $\phi:Y\rightarrow\overline{Y}$ is a composition of
type III birational contractions whose exceptional divisors
$E_1,\ldots,E_k$ are contracted down to a unique smooth
irreducible curve $C$ of compound Du Val singularities of type
$cA_{k}$. Examples given in \cite{KMP96} are 3--dimensional
hypersurfaces or complete intersections in weighted projective
spaces where birational contractions $\phi$'s are induced by
morphisms globally defined between the weighted projective spaces.
For each example a \emph{non--toric} deformation direction for $Y$
is exhibited, producing a deformation $\phi
':Y'\rightarrow\overline{Y}'$ of $\phi$ which is now a
\emph{small} birational contraction (a composition of type I
contractions). Moreover $\Sing ( \overline{Y}')$ turns out to be
composed only by nodes. Then $T$ deforms to a conifold transition
$T'(Y',\overline{Y}',\widetilde{Y})$ as follows:
\begin{equation}\label{T'}
    \xymatrix{Y\ar@/^3pc/
    @{.>}[drr]^T\ar[r]_{\phi}\ar@{<~>}[dd]^{\text{non--toric}}&
              \overline{Y}\ar@{<~>}[dr]\ar@{<~>}[dd]&\\
               &&\widetilde{Y}\\
              Y'\ar@/_3pc/ @{.>}[urr]^{T'}_{\text{conifold}}\ar[r]^{\phi '}&
              \overline{Y}'\ar@{<~>}[ur]&}
\end{equation}
In particular $E_1,\ldots,E_k$ are deformed to $\binom{k+1}{2}$
collections of $2g-2$ homologous rational curves in $Y'$, where
$g$ is the genus of $C$, and $C$ is deformed to
$N={\binom{k+1}{2}}(2g-2)$ nodes in $\overline{Y}'$. Since
$E_1,\ldots,E_k$ span a $k$--dimensional subspace of $H_4(Y)$, the
$N$ rational curves in $Y'$ span a $k$--dimensional subspace of
$H_2(Y')$. Setting $c=N-k$ one can then recover parameters $N,k,c$
for the given non--conifold transition $T$. If $T'$ admits a
reverse conifold transition $T'^{\circ}$ (of parameters $N,c,k$),
as Conjecture \ref{cong.reverse-conifold} predicts, then the
latter admits also $T$ as reverse transition. Therefore:
\begin{itemize}
    \item \emph{it can happen that a conifold transition of parameters
    $N,k,c$ admits a non--conifold reverse transition whose
    birational morphism contracts $c$ exceptional divisors down
    to a smooth irreducible curve of genus}
    \begin{equation*}
    g=1+\frac{N}{2\binom{c+1}{2}}\ ;
    \end{equation*}
\end{itemize}
This fact do not contradicts Conjecture
    \ref{cong.reverse-conifold} if the following one is true:
\begin{conjecture}
A geometric transition $T(Y,\overline{Y},\widetilde{Y})$
satisfying some good condition \emph{(e.g. such that $\phi$
contracts $k$ exceptional divisors down to a smooth curve of genus
$g>1$ whose points are $cA_k$ singularities)} can be deformed to a
conifold transition $T'(Y',\overline{Y}',\widetilde{Y})$ like in
diagram (\ref{T'}).
\end{conjecture}

\subsection{Toric degenerations: conifold transitions to construct mirror manifolds}
Methods in \cite{BC-FKvS98} were generalized in \cite{BC-FKvS00}
to complete intersections in partial flag manifolds giving a
conjectural approach to produce examples verifying Conjecture
\ref{cong.reverse-conifold}. On the other hand their method
describes a conjectural procedure \emph{to generalize the mirror
construction for \cy complete intersections in toric Fano
varieties, given in \cite{Batyrev94}, \cite{Batyrev-vStraten95}
and \cite{Borisov}, to the case of \cy complete intersections in
non--toric Fano varieties}. A further generalization of this
construction is given in \cite{Batyrev04}. Main ideas are the
following.

\begin{definition}[\cite{Batyrev04}, Definition 3.1]
Let $X\subset\P ^m$ be a smooth Fano variety of dimension $n$. A
normal Gorenstein toric Fano $P\subset\P ^m$ is called a
\emph{small toric degeneration of $X$}, if there exists a Zariski
open neighborhood $U$ of $0\in\C$ and an irreducible subvariety
$\mathcal{X}\subset\P ^m\times U$ such that the morphism
$\pi:\mathcal{X}\rightarrow U$ is flat and the following
conditions hold:
\begin{enumerate}
    \item the fiber $X_t:= \pi ^{-1}(t)\subset\P ^m$ is smooth for
    all $t\in U\setminus \{ 0\}$;
    \item the special fibre $X_0:=\pi ^{-1}(0)\subset\P ^m$ has at
    worst Gorenstein terminal singularities and $X_0\cong P$;
    \item the canonical homomorphism $\Pic
    (\mathcal{X}/U)\rightarrow\Pic(X_t)$ is an isomorphism for all
    $t\in U$.
\end{enumerate}
\end{definition}

\begin{examples}
\begin{enumerate}
    \item In \cite{BC-FKvS98} it is shown that the Grassmannian
    $X:=\mathbb{G}(r,s)$, embedded in $\P ^{\binom{s}{r}-1}$ by the
    usual Pl\"{u}cker embedding, admits a small toric degeneration
    $P:=P(r,s)\subset \P ^{\binom{s}{r}-1}$.
    \item In \cite{BC-FKvS00} it is proved that the partial flag
    manifold $X:=F(n_1,\ldots,n_k,n)$ with its Pl\"{u}cker
    embedding in $\P ^m$ admits a small toric degeneration $P\subset \P
    ^m$.
    \item In \cite{Batyrev04} the toric hypersurface $P$, given by the
    following homogeneous equation of degree $d$ in $\P^n$
    \begin{equation*}
    z_1\cdots z_d = z_{d+1}\cdots z_{2d}
    \end{equation*}
    where $n\geq 2d-2$, is proved to be a small toric degeneration of the
    generic smooth Fano hypersurface $X$ of degree $d$ in $\P^n$.
\end{enumerate}
\end{examples}

\begin{remark}
For all the previous examples, $\Sing P$ has codimension at least
3. Moreover the codimension 3 part of $\Sing P$ consists of
ordinary double points.
\end{remark}

Let now $H$ be a generic complete intersection in $\P^m$ cutting
on a smooth Fano variety $X\subset \P^m$ a smooth \cy variety $Y$.
If $X$ admits a small toric degeneration $P\subset\P^m$ and
$\overline{Y}:=H\cap P$ then $\Sing \overline{Y}$ has codimension
at least 3. In particular \emph{if $\dim Y = 3 =\dim\overline{Y}$
then $\Sing \overline{Y}$ consists only of isolated nodes}. Let
$\widehat{P}$ be a simultaneous desingularization of $P$ given by
a suitable subdivision of the fan associated with $P$. Then the
birational morphism $\widehat{P}\rightarrow P$ induces a
desingularization $\widehat{Y}\rightarrow Y$. We have then a
geometric transition $T(\widehat{Y},\overline{Y},Y)$ which is
conifold when $\dim Y=3$.
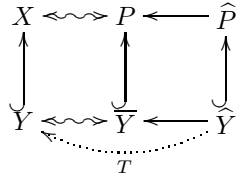
\begin{equation*}
    \xymatrix{X\ar@{<~>}[r]&P&\widehat{P}\ar[l]\\
              Y\ar@{^{(}->}[u]\ar@{<~>}[r]&\overline{Y}\ar@{^{(}->}[u]&
              \widehat{Y}\ar[l]\ar@/^1pc/ @{.>}[ll]^T\ar@{^{(}->}[u]}
\end{equation*}
The mirror partner of $\widehat{Y}$ given by the construction of
\cite{Batyrev-vStraten95} and \cite{Borisov}, is a complete
intersection $\widehat{Y}^{\circ}$ in the dual Fano toric variety
$\widehat{P}^{\circ}$ obtained by polarity on associate polytopes.
The main point is that the embedding $\Pic P \hookrightarrow \Pic
\widehat{P}$ suggests, via monomial--divisor correspondence
\cite{AGM93}, a canonical way to specialize $\widehat{Y}^{\circ}$
to a singular $\overline{Y}^{\circ}$. Let $Y^{\circ}\rightarrow
\overline{Y}^{\circ}$ be a minimal desingularization. The
situation is then the following
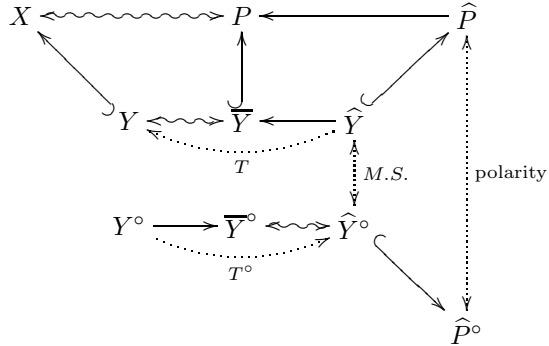
\begin{equation*}
    \xymatrix{X\ar@{<~>}[rr]&&P&&\widehat{P}\ar[ll]\\
              &Y\ar@{^{(}->}[lu]\ar@{<~>}[r]&\overline{Y}\ar@{^{(}->}[u]&
              \widehat{Y}\ar[l]\ar@/^1pc/ @{.>}[ll]^T\ar@{^{(}->}[ur]\ar@{<.>}[d]^{M.S.}&\\
              &Y^{\circ}\ar@/_1pc/
              @{.>}[rr]_{T^{\circ}}\ar[r]&\overline{Y}^{\circ}\ar@{<~>}[r]&\widehat{Y}^{\circ}\ar@{^{(}->}[dr]&\\
              &&&&\widehat{P}^{\circ}\ar@{<.>}[uuu]_{\text{polarity}}}
\end{equation*}
and \emph{$Y^{\circ}$ is conjectured to be a mirror partner of $Y$
and $T^{\circ}$ be a reverse transition of $T$}. In particular for
all the given 3--dimensional examples verifying this conjecture
(see \cite{BC-FKvS98}) $T^{\circ}$ turns out to be a conifold
transition like $T$.

\subsection{Mirror partners of rigid \cy 3--folds via geometric transitions}

Let $Y$ be a rigid \cy 3--fold as in Remark \ref{rigid-cy}. At
least from the mathematical point of view, the reverse transition
Conjecture \ref{cong.reverse} gives an answer to which non--\ka
\cy 3-fold $Y^{\circ}$ should be a mirror partner of $Y$. In fact
\begin{itemize}
    \item if there exists a geometric transition
    $T(Y,\overline{Y},\widetilde{Y})$ then
    $h^{2,1}(\widetilde{Y})>h^{2,1}(Y)=0$, since $\widetilde{Y}$
    cannot be rigid,
    \item let $\widetilde{Y}^{\circ}$ be a mirror partner of
    $\widetilde{Y}$ and
    $T^{\circ}(\widetilde{Y}^{\circ},\overline{Y}^{\circ},Y^{\circ})$
    be a reverse transition of $T$,
    \item then $Y^{\circ}$ should be a mirror partner of $Y$ like
    in diagram (\ref{reverse}).
\end{itemize}
If $T$ and $T^{\circ}$ are both conifold then, from the physical
point of view, the previous procedure suggests that the mirror
dual of a string theory compactified to a rigid \cy 3--fold can be
obtained by a suitable composition of black hole condensations and
mirror symmetry (over non--rigid \cy 3--folds).

\section{Further physical dualities and transitions}

The local conifold transition
\begin{equation}\label{conifold locale}
    \xymatrix@1{\mathcal{O}_{\P^1}(-1)\oplus\mathcal{O}_{\P^1}(-1)\ar@/_1pc/ @{.>}[rr]_T
              \ar[r]&\overline{U}\ar@{<~>}[r]& T^*S^3}
\end{equation}
studied in section \ref{analisi locale}, has been recently
considered as the geometric set up of a new conjectured
open/closed string duality.

\noindent More precisely, at the beginning, in 1974, G.~t'Hooft
conjectured that large $N$ gauge theories are dual to closed
string theories, \cite{t'Hooft74}. Later, in 1992, E.~Witten showed
that a particular kind of gauge theory, namely a $\su (N)$ (or
$U(N)$) Chern--Simons gauge theory on the 3--sphere $S^3$, is
equivalent to an open string theory on $T^* S^3$ with D--branes
wrapped on $S^3$, \cite{Witten95}. In 1998, R.~Gopakumar and C.~Vafa
conjectured that, for large $N$, a $\su (N)$ ($U(N)$)
Chern--Simons gauge theory is dual to a closed string theory
``compactified" to the local \cy 3--fold
$\mathcal{O}_{\P^1}(-1)\oplus\mathcal{O}_{\P^1}(-1)$,
\cite{Gopakumar-Vafa99} and \cite{Ooguri-Vafa00}. Composing all
these dualities gives an open/closed string duality modelled on
the local conifold transition (\ref{conifold locale}). For all the
details, the interested reader is referred to original papers, and
to \cite{GR02} for a survey on these topics and more references.

\noindent The concept of reverse transition, introduced in the
previous section, applied to such an open/closed string duality,
suggests a further duality on the mirror theories. This was
proposed in \cite{Aganagic-Vafa01}.

\noindent Examples of similar dualities, geometrically realized by
less elementary conifold transitions than (\ref{conifold locale}),
are given in \cite{DFG03}. A reverse transition of one of them is
described in the recent paper \cite{Forbes04}.

\end{document}